\documentclass[english,12pt,onecolumn]{IEEEtran}
\usepackage[T1]{fontenc}
\usepackage[latin9]{inputenc}
\usepackage{color}
\usepackage{babel}
\usepackage{amsmath}
\usepackage{amsthm}
\usepackage{amssymb}
\usepackage{graphicx}
\usepackage{setspace}
\usepackage[unicode=true,pdfusetitle,
 bookmarks=true,bookmarksnumbered=true,bookmarksopen=false,
 breaklinks=true,pdfborder={0 0 1},backref=false,colorlinks=true]
 {hyperref}
\hypersetup{
 pdfborderstyle=}

\makeatletter
\numberwithin{figure}{section}
\theoremstyle{plain}
\newtheorem{thm}{\protect\theoremname}
\theoremstyle{definition}
\newtheorem{problem}{\protect\problemname}
\theoremstyle{plain}
\newtheorem{cor}{\protect\corollaryname}
\theoremstyle{plain}
\newtheorem{prop}{\protect\propositionname}
\theoremstyle{plain}
\newtheorem{lem}{\protect\lemmaname}
\theoremstyle{remark}

\makeatother

\providecommand{\corollaryname}{Corollary}
\providecommand{\lemmaname}{Lemma}
\providecommand{\problemname}{Problem}
\providecommand{\propositionname}{Proposition}
\providecommand{\theoremname}{Theorem}

\begin{document}
\title{On Average Distance, Level-1 Fourier Weight, and Chang's Lemma}
\author{Lei Yu\thanks{L. Yu is with the School of Statistics and Data Science, LPMC, KLMDASR,
and LEBPS, Nankai University, Tianjin 300071, China (e-mail: leiyu@nankai.edu.cn).
}
}
\maketitle
\begin{abstract}
We study the maximum level-$1$ Fourier weight of Boolean functions,
which is equivalent to the minimum average-distance problem on the
hypercube. We first determine the dimension-free optimum at sufficiently
small densities: there exists a universal $a_{0}>0$ such that
 the maximum  level-$1$ Fourier weight for Boolean functions of mean smaller than $a_{0}$   is asymptotically attained by Hamming balls as the dimension $n\to \infty$. The key ingredient is an eventual Gaussian
stop-loss domination inequality for normalized Rademacher sums. We then
use an induction argument to improve the classical level-$1$ bound
(Chang's lemma) in both the small- and large-density regimes. We apply
these estimates to strengthen the Friedgut--Kalai--Naor theorem, study
the corresponding average-distance problem in Euclidean space, and
derive a sharp form of Chang's original lemma for $\mathbb{F}_{2}^{n}$:
Hamming balls maximize the dimension of the span of the large Fourier
coefficients.
\end{abstract}


\section{\label{sec:Introduction}Introduction}

Fourier analysis of Boolean functions plays an important role in coding
theory, theoretical computer science, and mathematics. In this paper,
we focus on maximizing the level-$1$ Fourier weight of Boolean functions,
which is equivalent to an anticoding problem: minimizing the
average pairwise distance of a code. The problem is formulated
as follows. 

Consider the Fourier basis $\left\{ \chi_{S}\right\} _{S\subseteq[n]}$
on the hypercube $\{\pm1\}^{n}$ with $\chi_{S}(\mathbf{x}):=\prod_{i\in S}x_{i}$
for $\mathbf{x}:=(x_{1},x_{2},\ldots,x_{n})\in\{\pm1\}^{n}$ and $S\subseteq[n]:=\{1,2,\ldots,n\}$.
Then for a Boolean function $f:\{\pm1\}^{n}\to\{0,1\}$, its Fourier
coefficients are 
\begin{equation}
\hat{f}_{S}:=\mathbb{E}[f(\mathbf{X})\chi_{S}(\mathbf{X})],\;S\subseteq[n],\label{eq:-33}
\end{equation}
where $\mathbf{X}\sim\mu^{\otimes n}:=\mathrm{Unif}\{\pm1\}^{n}$.
The Fourier expansion of a Boolean function $f$ (cf. \cite[Equation (1.6)]{ODonnell14analysisof})
is 
\begin{align*}
 & f(\mathbf{x})=\sum_{S\subseteq[n]}\hat{f}_{S}\chi_{S}(\mathbf{x}).
\end{align*}
The \emph{level-$k$ Fourier weight} of $f$ is defined as 
\begin{align*}
\mathbf{W}_{k}[f] & :=\sum_{S:|S|=k}\hat{f}_{S}^{2},\quad k\in[n],
\end{align*}
which is the energy of the Fourier coefficients on the Hamming sphere
of radius $k$. For a Boolean function, we also denote $\mathbf{W}_{k}[f]$
by $\mathbf{W}_{k}[\operatorname{supp}(f)]$, or shortly, $\mathbf{W}_{k}$. By
definition and Parseval's theorem, it is easily seen that 
\[
\mathbf{W}_{0}=a^{2}\quad\textrm{and}\quad\sum_{k=0}^{n}\mathbf{W}_{k}=a,
\]
where $a=\mathbb{E}[f]$. The most important Fourier weight besides
$\mathbf{W}_{0}$ is the degree-$1$ Fourier weight 
\begin{align*}
\mathbf{W}_{1} & :=\sum_{i=1}^{n}\hat{f}_{i}^{2},
\end{align*}
where for short, we denote $\hat{f}_{i}=\hat{f}_{\{i\}}$ for $i\in[n]$.
Estimating Fourier coefficients and Fourier weights of Boolean
functions is a central topic with
many applications in theoretical computer science, Fourier analysis,
coding theory, and combinatorics \cite{ODonnell14analysisof,green2008boolean,chang2002polynomial,friedgut2002boolean,defant2019fourier}.
A classical upper bound on $\mathbf{W}_{1}$ is the following
level-$1$ inequality \cite{ODonnell14analysisof}, which is also known
as Chang's lemma. Let
\[
\mathcal D_n:=\{0,2^{-n},2\cdot2^{-n},\ldots,1\},
\qquad
\mathcal D:=\bigcup_{n\ge0}\mathcal D_n
\]
be the admissible densities in dimension $n$ and the set of dyadic
densities, respectively. For $a\in\mathcal D_n$, define
\begin{align}
W^{(n)}(a) & :=\max_{\textrm{Bool }f:\mathbb{E}f=a}\mathbf{W}_{1}[f],\label{eq:Wn}\\
W(a) & :=\sup_{n:\,a\in\mathcal D_n}W^{(n)}(a),\qquad a\in\mathcal D.\label{eq:W}
\end{align}
Section~\ref{sec:small-density} shows that $W$ has a unique continuous
extension to $[0,1]$; throughout the paper we use the same symbol for
this extension. By complement symmetry, $W(1-a)=W(a)$, so it suffices
to consider $a\in[0,1/2]$.
\begin{thm}[Level-1 Weight Bound]
\label{thm:level-1-bound}
For every $a\in(0,1/2]$,
\[
W(a)\le2a^{2}\ln\frac{1}{a}.
\]
\end{thm}
This bound follows readily from hypercontractivity or from the entropy
method. This inequality has many applications in the analysis
of Boolean functions and additive combinatorics.  As the dimension
$n\to\infty$, the degree-$1$ Fourier weight $\mathbf{W}_{1}$ of
Hamming balls of density $a$ converges to $J(a)=I^{2}(a)$, where $I(a)=\phi(\Phi^{-1}(a))$
is the Gaussian isoperimetric profile with $\phi$ and $\Phi$ respectively
denoting the standard Gaussian density and cumulative distribution
function (cdf). Since $J(a)\sim2a^{2}\ln\frac{1}{a}$ as $a\to0$,
the lemma above implies that, as $a\to0$, indicators of Hamming balls
asymptotically maximize $\mathbf{W}_{1}$ over all Boolean functions
of the same size. Since finite-dimensional Hamming balls approach the
value $J(a)$, it is natural to ask whether this Gaussian/Hamming-ball
limit is the exact dimension-free optimum for all sufficiently small
fixed densities. 
\begin{problem}
\label{prob:W}Is there a number $a_{0}>0$ such that $W(a)=J(a)$
for every $a\in(0,a_{0}]$?
\end{problem}
It is known that for $a\in\{1/2,1/4\}$, subcubes maximize
$\mathbf{W}_{1}$ exactly, so the Gaussian/Hamming-ball value is not
optimal at these larger densities. We answer Problem~\ref{prob:W}
affirmatively in Section~\ref{sec:small-density}: there exists a
universal $a_{0}>0$ such that $W(a)=J(a)$ for every $0<a\le a_{0}$.
Moreover, the supremum is not attained by any finite-dimensional
Boolean function in this regime. We then develop quantitative upper
bounds that remain useful beyond the range covered by this exact
small-density theorem.

In fact, if we define the average distance $\mathbf{D}(A)$ of $A\subseteq\{\pm1\}^{n}$
as the average of the Hamming distance of every pair of elements in
$A$, then the following relation \cite{shutao1998average} holds:
\begin{align}
\mathbf{W}_{1}[A] & =a^{2}(n-2\mathbf{D}(A)).\label{eq:WD}
\end{align}
where $a=\mu^{\otimes n}(A)$. In other words, given the size of $A$,
maximizing the level-1 Fourier weight of $1_{A}$ is equivalent to
minimizing the average distance of $A$. The latter problem was posed
by Ahlswede and Katona \cite{ahlswede1977contributions}, and it is
an ``anticoding problem.'' Ahlswede and Althöfer \cite{ahlswede1994asymptotic}
considered the case in which the size of $A$ increases exponentially
in $n$ and the exponent is strictly between $0$ and $1$. They showed
that Hamming spheres asymptotically minimize the average distance
as $n\to\infty$. Using a linear programming approach, Mounits \cite{mounits2008lower}
studied sets whose sizes are linear in $n$ (i.e., sets with ``small''
sizes). He showed that when the set size is $2n$, the asymptotic
value of the minimum average distance is $\frac{5}{2}$ as $n\to\infty$.
Althöfer and Sillke \cite{althofer1992average}, Fu, Xia, together
with other authors \cite{shutao1998average,fu1997expectation,fu1999hamming,fu2001minimum},
as well as Mounits \cite{mounits2008lower}, proved various bounds
on the minimum average distance, which are sharp in certain regimes
when the code size is ``large'' (e.g., the size is $2^{n-1}$ or
$2^{n-2}$). In contrast, Theorem~\ref{thm:optimality} does not solve
the problem itself, but characterizes the structure of the optimal solutions.
That is, the optimal solutions to Ahlswede--Katona's problem (or
the maximum level-$1$ Fourier weight problem) must be linear-threshold
sets or half-spaces (i.e., the supports of LTFs). Thus, to solve Ahlswede--Katona's problem, it suffices to consider
linear-threshold sets. For large $a$, e.g., $a=1/2$ or $1/4$, subcubes
were shown to be exactly optimal \cite{althofer1992average,shutao1998average},
while for small $a$, e.g., exponentially small $a$, Hamming balls
were shown to be asymptotically optimal \cite{ahlswede1994asymptotic}.

\subsection{Our Contributions }
\begin{enumerate}
\item We resolve Problem~\ref{prob:W}. There exists a universal
$a_{0}>0$ such that
\[
W(a)=J(a)=I(a)^2,\qquad 0<a\le a_{0}.
\]
The main probabilistic ingredient is an eventual Gaussian stop-loss
comparison for normalized Rademacher sums. In fact, every
finite-dimensional Boolean function of mean $a\in(0,a_0]$ satisfies the
strict inequality $\mathbf{W}_{1}[f]<J(a)$, so the dimension-free value
is approached only in the limit of growing dimension.
\item We prove two new upper bounds on the level-$1$ Fourier weight. Our
first bound is $W(a)\le\chi(a)$ for $a\in[0,1/2]$, where
\begin{equation}
\chi(a)=\begin{cases}
I_{w}^{2}(a), & a\in[0,T]\\
a/2, & a\in(T,1/2]
\end{cases},\label{eq:J-1}
\end{equation}
with $T=0.21$, $I_{w}(a):=wI\left(\frac{a}{w}\right)$, and $w$
denoting the unique solution to $I_{w}(T)=\sqrt{T/2}$. This bound
improves Chang's bound in Theorem~\ref{thm:level-1-bound}, as well as
previous bounds for $a\in(0,T)$. From it we derive a second bound that
improves existing estimates for $a\ge1/4$ and is numerically
asymptotically tight as $a\to1/2$. Through \eqref{eq:WD}, these estimates
also yield bounds for the minimum average-distance problem.
\item We apply the new upper bound to improve the
Friedgut--Kalai--Naor theorem for balanced Boolean functions. If
$f$ is balanced and $\max_{i\in[n]}|\hat f_i|=\beta$, then
$\mathbf{W}_{1}[f]\le\beta^{2}+\chi(\frac12-\beta)$. This estimate is
asymptotically sharp in a suitable sense as $\beta\to1/2$.
\item We study the analogue of the average-distance problem in Euclidean
space and show that Euclidean balls minimize average distance among all
sets of a prescribed measure.
\item We derive a sharp form of Chang's original lemma for
$\mathbb{F}_{2}^{n}$: Hamming balls maximize the dimension of the space
spanned by large Fourier coefficients.
\end{enumerate}

\subsection{Preliminaries}

It is known that every maximizer of $\mathbf{W}_{1}$ must be a linear-threshold function. For a Boolean function $f$, denote
\begin{align*}
f_{\ge}(\mathbf{x}) & =1\Big\{\mathbf{x}:\sum_{i=1}^{n}\hat{f}_{i}x_{i}\ge b\Big\}\\
f_{>}(\mathbf{x}) & =1\Big\{\mathbf{x}:\sum_{i=1}^{n}\hat{f}_{i}x_{i}>b\Big\},
\end{align*}
for a suitable $b$ such that $\mathbb{E}f_{>}\le a\le\mathbb{E}f_{\ge}$.
 These are linear-threshold functions. 
\begin{thm}[Self-Consistency]
\label{thm:optimality}Any Boolean function $f:\{\pm1\}^{n}\to\{0,1\}$
maximizing $\mathbf{W}_{1}[f]$ over all Boolean functions such that $\mathbb{E}f=a$
satisfies 
\begin{equation}
f_{>}\le f\le f_{\ge}.\label{eq:-34-1-1}
\end{equation}
Moreover, the set of optimal Boolean functions is $\{\textrm{Bool }\ell:f_{>}\le\ell\le f_{\ge},\mathbb{E}\ell=a\}$. 
\end{thm}
This theorem was proven independently in \cite{yu2023entropy} and
\cite{heilman2025fourier}. A short proof is given here, and an
alternative proof is provided in Appendix~\ref{sec:-Proof-of}. 
\begin{IEEEproof}
Let $A$ be the support of $f$. Observe that $\hat{f}_{i}=\mathbb{E}[f(\mathbf{X})X_{i}]=\frac{a}{|A|}\sum_{\mathbf{x}\in A}x_{i}$.
So, $a^{-1}(\hat{f}_{i})_{i\in[n]}$ is the barycenter of $A$, and
$a^{-1}\sqrt{\mathbf{W}_{1}[f]}$ is the Euclidean distance from the
origin to this barycenter. Given $a$, maximizing $\mathbf{W}_{1}[f]$
over all Boolean $f$ of mean $a$ is equivalent to maximizing the
distance from the origin to the barycenter of a set $A$ over all
$A$ of measure $a$. Obviously, any optimal set $A$ with barycenter
$\mathbf{w}$ for the latter problem must satisfy $\{\left\langle \mathbf{w},\mathbf{x}\right\rangle >b\}\subseteq A\subseteq\{\left\langle \mathbf{w},\mathbf{x}\right\rangle \ge b\}$
for some $b$ chosen properly to ensure the existence of such $A$
satisfying $\mu(A)=a$. Otherwise, one can move points
in $A$ from $\{\left\langle \mathbf{w},\mathbf{x}\right\rangle <b\}$
to $\{\left\langle \mathbf{w},\mathbf{x}\right\rangle \ge b\}$, which
produces a larger distance from the origin to the new barycenter. 
\end{IEEEproof}
Thus, it suffices to consider only linear-threshold functions
when maximizing $\mathbf{W}_{1}$. The indicators of subcubes and Hamming
balls are both linear threshold functions. As a consequence of Theorem
\ref{thm:optimality}, both $W^{(n)}(a)$ and $W(a)$ are nondecreasing
in $a\le1/2$. 
\begin{cor}
Given $n$, $W^{(n)}(a)$ is strictly increasing in $a\le1/2$. Consequently, $W(a)$ is nondecreasing  in $a\le1/2$. 
\end{cor}

\begin{IEEEproof}
Let $a<1/2$. Suppose that $f$ attains $W^{(n)}(a)$. Let $A$ be
the support of $f$. Without loss of generality, we assume $\hat{f}_{i}\ge0$
for all $i\in[n]$, since otherwise, we can change the sign of the
$i$-th coordinate for all points in $A$. So, it satisfies \eqref{eq:-34-1-1}
with $b\ge0$. Let $\mathbf{y}$ be a point in $A^{c}$ having the
largest $\sum_{i=1}^{n}\hat{f}_{i}y_{i}$. So, $\sum_{i=1}^{n}\hat{f}_{i}y_{i}\ge0$.
Moreover, $g(\mathbf{x})=f(\mathbf{x})+1\{\mathbf{x}=\mathbf{y}\}$
has Fourier coefficients $\hat{g}_{i}=\hat{f}_{i}+2^{-n}y_{i}$. Hence,
\[
\mathbf{W}_{1}[g]=\mathbf{W}_{1}[f]+n4^{-n}+2^{1-n}\sum_{i=1}^{n}\hat{f}_{i}y_{i}>\mathbf{W}_{1}[f],
\]
which implies $W^{(n)}(a+2^{-n})>W^{(n)}(a)$ for any $a<1/2$. 
Hence, $W^{(n)}(a)$ is strictly increasing in $a\le1/2$. 
\end{IEEEproof}
As mentioned above, Ahlswede and Althöfer \cite{ahlswede1994asymptotic}
proved that for the case that $\mu^{\otimes n}(A)$ vanishes exponentially
in $n$, Hamming balls or spheres asymptotically minimize the average
distance as $n\to\infty$. By the equivalence in \eqref{eq:WD} and
by tensorization, their result implies the following bound
for finite dimension. 
\begin{thm}[Strong Level-1 Weight Bound]
\label{thm:expW1} Let $f:\{\pm1\}^{n}\to\{0,1\}$ have expectation
$a$. Then, 
\begin{equation}
\mathbf{W}_{1}[f]\le na^{2}\Big(1-2H^{-1}(1-\frac{1}{n}\log_{2}\frac{1}{a})\Big)^{2},\label{eq:-32-1}
\end{equation}
where $H^{-1}$ is the inverse of the binary entropy function $H:t\in[0,1]\mapsto-t\log_{2}t-(1-t)\log_{2}(1-t)$.
This bound is asymptotically sharp, since a sequence of Hamming balls
$B_{r}^{(n)}=\{\mathbf{x}:\frac{1}{n}\sum_{i=1}^{n}x_{i}\le r\}$
for $r$ properly chosen such that $-\frac{1}{n}\log_{2}\mu^{\otimes n}(B_{r}^{(n)})\to\alpha$
as $n\to\infty$ satisfies $\mathbf{W}_{1}[B_{r}^{(n)}]\sim na^{2}\Big(1-2H^{-1}(1-\frac{1}{n}\log_{2}\frac{1}{a})\Big)^{2},$
where $a=2^{-n\alpha}$. 
\end{thm}
Here we provide a direct proof for this result. 
\begin{IEEEproof}
Let $\lambda_{\mathbf{X}}:=\mu^{\otimes n}(\cdot|A)$, where $A$
is the support of $f$. Then, $\hat{f}_{i}=a(\lambda_{X_{i}}(1)-\lambda_{X_{i}}(-1)).$
Denote 
\begin{equation}
\beta:=\mathbf{W}_{1}[f]\big/a^{2}=\sum_{i=1}^{n}\|\lambda_{X_{i}}-\mu\|_{1}^{2}.\label{eq:-7}
\end{equation}
Denote the relative entropy as $D(\lambda\|\mu)=\sum_{\mathbf{x}}\lambda(\mathbf{x})\log_{2}\frac{\lambda(\mathbf{x})}{\mu(\mathbf{x})}$
for two distributions $\lambda,\mu$. Denote for $b\in[0,1]$, 
\begin{equation}
\Upsilon(b):=\inf_{\lambda:\|\lambda-\mu\|_{1}=b}D(\lambda\|\mu)=D(\frac{1-b}{2}\|\frac{1}{2}),\label{eq:upsilon}
\end{equation}
which is convex in $b$. Observe that 
\begin{align*}
\log\frac{1}{a} & =D(\lambda_{\mathbf{X}}\|\mu^{\otimes n})\ge\inf_{(\lambda_{X_{i}}):\sum_{i=1}^{n}\|\lambda_{X_{i}}-\mu\|_{1}^{2}=\beta}\sum_{i=1}^{n}D(\lambda_{X_{i}}\|\mu)\\
 & =\inf_{b_{i}\ge0:\sum_{i=1}^{n}b_{i}=\beta}\sum_{i=1}^{n}\Upsilon(\sqrt{b_{i}})=n\Upsilon(\sqrt{\beta/n})=nD\Big(\frac{1-\sqrt{\beta/n}}{2}\|\frac{1}{2}\Big),
\end{align*}
where the first inequality uses the superadditivity of relative entropy.
This is exactly \eqref{eq:-32-1}. 

The asymptotic sharpness of \eqref{eq:-32-1} can be proven by using
the large deviations theorem to estimate the exponent of $\mu^{\otimes n}(B_{r}^{(n)})$,
using the Gibbs conditioning principle to estimate $\lambda_{X_{i}}$,
and plugging the estimation of $\lambda_{X_{i}}$ into the relation
in \eqref{eq:-7}. 
\end{IEEEproof}

\subsection{Organization}

The rest of the paper is organized as follows.
Section~\ref{sec:small-density} proves the exact identity $W(a)=J(a)$
for all sufficiently small densities, via an eventual Gaussian
stop-loss comparison for normalized Rademacher sums.
Section~\ref{sec:improved-bounds} develops improved upper bounds on the
level-$1$ Fourier weight in the regimes $a<1/4$ and $a>1/4$.
Section~\ref{sec:fkn} applies these bounds to the
Friedgut--Kalai--Naor theorem. Section~\ref{sec:euclidean} studies the
Euclidean average-distance problem, and Section~\ref{sec:chang} treats
Chang's original lemma. Section~\ref{sec:open-problems} records two
remaining questions. Auxiliary proofs are deferred to the appendices.

\section{Small-Density Gaussian Optimality}\label{sec:small-density}

We now resolve Problem~\ref{prob:W}. The key step is an eventual
Gaussian stop-loss comparison for normalized Rademacher sums; a simple
Fourier duality argument then converts this probabilistic comparison
into the exact value of the dimension-free level-$1$ Fourier weight at
sufficiently small densities. Recall that $t\mapsto\mathbb{E}(X-t)_+$
is the stop-loss transform of $X$; see, e.g.,
\cite[Sec.~4.A]{ShakedShanthikumar2007} for its role in stochastic orders.

We retain the notation introduced above: $\phi$ and $\Phi$ denote the
standard Gaussian density and cumulative distribution function,
respectively, $I(a)=\phi(\Phi^{-1}(a))$, and $J(a)=I(a)^2$. We also write
\[
\overline{\Phi}(x):=1-\Phi(x),
\qquad
u_+:=\max\{u,0\}.
\]
Let $Z\sim N(0,1)$. All Rademacher variables $X_i$ below are independent and
take the values $\pm1$ with equal probability.

\begin{thm}[Eventual Gaussian Stop-Loss Domination]
\label{thm:stoploss}
There exists a universal constant $T_0<\infty$ such that, for every
$n\ge1$, every $c_1,\ldots,c_n\in\mathbb{R}$ satisfying
$\sum_{i=1}^n c_i^2=1$, and every $t\ge T_0$,
\begin{equation}\label{eq:sd-stoploss}
\mathbb{E}\left(\sum_{i=1}^n c_i X_i-t\right)_+
<\mathbb{E}(Z-t)_+.
\end{equation}
\end{thm}

\begin{thm}[Small-Density Gaussian Optimality]
\label{thm:small-density}
Let $T_0$ be as in Theorem~\ref{thm:stoploss}, and set
$a_0=\overline{\Phi}(T_0)>0$. Then
\begin{equation}\label{eq:sd-main}
W(a)=J(a),\qquad 0<a\le a_0.
\end{equation}
Moreover, every finite-dimensional Boolean function $f$ of mean
$a\in(0,a_0]$ satisfies the strict inequality
\[
\mathbf{W}_{1}[f]<J(a).
\]
Thus $J(a)$ is a supremum over dimensions, rather than a
finite-dimensional maximum, throughout this interval.
\end{thm}

We next provide the proofs for these two theorems. 

\subsection{Normal approximation with a single kink}

Fix a universal Berry--Esseen constant $C_{\mathrm{BE}}$ such that independent
centered random variables $X_j$, with $s^2=\sum_j\mathbb{E} X_j^2>0$, satisfy
\begin{equation}\label{eq:sd-BE}
\sup_x\left|
\mathbb{P}\left(\frac{\sum_jX_j}{s}\le x\right)-\Phi(x)
\right|
\le \frac{C_{\mathrm{BE}}}{s^3}\sum_j\mathbb{E}|X_j|^3.
\end{equation}
We use only this standard form of the Berry--Esseen theorem; see, e.g.,
\cite[Ch.~3]{ChenGoldsteinShao2011}.

\begin{lem}[Signed-third-moment estimate]
\label{lem:normal-kink}
Let $Y_1,\ldots,Y_n$ be independent centered random variables such that
\[
\sum_i\mathbb{E} Y_i^2=1,
\qquad
|Y_i|\le r_i\le\frac12,
\qquad
\sum_i r_i^2\le B.
\]
Let $Y=\sum_iY_i$. Suppose that $h:\mathbb{R}\to\mathbb{R}$ is continuous,
$L$-Lipschitz, and $C^2$ on each side of a point $k\in[-1,1]$, with
\[
|h''(x)|\le M,\qquad x\ne k.
\]
Then
\begin{equation}\label{eq:sd-normal}
\left|\mathbb{E} h(Y)-\mathbb{E} h(Z)\right|
\le C_B\left[
L\sum_i|\mathbb{E} Y_i^3|+(L+M)\sum_i r_i^4
\right],
\end{equation}
where one may take
\begin{equation}\label{eq:sd-CB}
C_B=100+12BC_{\mathrm{BE}}.
\end{equation}
\end{lem}

\begin{IEEEproof}
Let $F$ be the bounded solution of the normal Stein equation
\begin{equation}\label{eq:sd-stein}
F'(x)-xF(x)=h(x)-\mathbb{E} h(Z).
\end{equation}
The existence and uniqueness of the bounded solution, together with the
standard Stein-factor estimates used below, are classical; see
\cite[Sec.~2.1]{Ross2011} and \cite[Lem.~2.4]{ChenGoldsteinShao2011}.
The one-sided derivatives of $h$ at $k$ exist. 
For a function with one-sided limits at $k$, write
\[
[f]_k:=f(k+)-f(k-)
\]
for its jump at $k$.
Write
\[
J_h=h'(k+)-h'(k-)=[h']_k,\qquad |J_h|\le2L.
\]
We first establish the bounds
\begin{equation}\label{eq:sd-steinregularity}
\|F''\|_\infty\le4L,
\qquad
[F'']_k=J_h,
\qquad
\sup_{x\ne k}|F'''(x)|\le128(L+M).
\end{equation}
At $k$, either one-sided value may be assigned to $F''$.

\medskip
\noindent\emph{Step 1: Stein regularity.}
For any Lipschitz function $u$, its bounded Stein solution has the
Ornstein--Uhlenbeck representation
\[
F_u(x)=-\int_0^1
\mathbb{E}\,u'\!\left(sx+\sqrt{1-s^2}\,Z\right)\,ds.
\]
This semigroup representation, as well as the Gaussian integration-by-parts
identity used below, is standard; see \cite[Secs.~2.8--2.9 and~3.2]{NourdinPeccati2012}.
Gaussian integration by parts therefore gives
\begin{equation}\label{eq:sd-steinbasic}
\|F_u\|_\infty\le\operatorname{Lip}(u),
\qquad
\|F_u'\|_\infty
\le\mathbb{E}|Z|\,\operatorname{Lip}(u)\le\operatorname{Lip}(u).
\end{equation}
We also use the elementary estimate
\begin{equation}\label{eq:sd-boundedstein}
\|g'\|_\infty\le2\|q\|_\infty
\end{equation}
for any bounded, locally absolutely continuous solution of
$g'-xg=q$ with bounded and centered (w.r.t. Gaussian measure  $\gamma$) $q$.
This follows the following argument. Since $\phi'(x)=-x\phi(x)$, the Stein equation
$
g'(x)-xg(x)=q(x)
$
gives
$
(g\phi)'=q\phi.
$
Using $\int_{\mathbb R}q\,d\gamma=0$, we therefore obtain the two
representations
\[
g(x)
=\frac{1}{\phi(x)}\int_{-\infty}^x q(t)\phi(t)\,dt
=-\frac{1}{\phi(x)}\int_x^\infty q(t)\phi(t)\,dt.
\]
For $x>0$, the second representation and the Gaussian Mills bound
$\bar\Phi(x)\le \phi(x)/x$ give
\[
|xg(x)|
\le \|q\|_\infty
\frac{x\bar\Phi(x)}{\phi(x)}
\le \|q\|_\infty.
\]
For $x<0$, the first representation gives the same conclusion by
symmetry. Hence
$
\|xg\|_\infty\le \|q\|_\infty.
$
This equation
then implies \eqref{eq:sd-boundedstein}. Bounds involving derivatives
are interpreted almost everywhere when necessary.

Differentiating \eqref{eq:sd-stein} once gives, almost everywhere,
\[
(F')'-xF'=F+h'.
\]
Thus $g:=F'$ solves the Stein equation $g'-xg=q$ with
$q:=F+h'$. Since $\operatorname{Lip}(h)\le L$, we have
$\|h'\|_\infty\le L$ a.e., while \eqref{eq:sd-steinbasic} gives
$\|F\|_\infty\le L$. Hence
\[
\|q\|_\infty\le \|F\|_\infty+\|h'\|_\infty\le2L.
\]
Applying \eqref{eq:sd-boundedstein} therefore yields
\begin{equation}
\|F''\|_\infty
\le2\|q\|_\infty
\le4L. \label{eq:sd-steinexpansion}
\end{equation}
The identity
$F''=F+xF'+h'$ also shows that $[F'']_k=J_h$.

Decompose
\begin{equation}
h(x)=J_h \cdot (x-k)_++r(x). \label{eq:decomp}
\end{equation}
Then $r'$ is continuous, $\operatorname{Lip}(r)\le3L$, and $|r''|\le M$
almost everywhere. The Stein solution for $r$ satisfies
\[
(F_r'')'-xF_r''=2F_r'+r''.
\]
Since $F_r''$ is bounded, we may apply
\eqref{eq:sd-boundedstein} to
$
g:=F_r'',
\,
q:=2F_r'+r'',
$
yielding that 
\[
\|F_r'''\|_\infty
\le 2\|2F_r'+r''\|_\infty
\le 4\|F_r'\|_\infty+2\|r''\|_\infty.
\]
By \eqref{eq:sd-steinbasic} and $\operatorname{Lip}(r)\le3L$, we have 
$
\|F_r'\|_\infty\le3L,
$
while $\|r''\|_\infty\le M$. Consequently,
\begin{equation}\label{eq:sd-smooththird}
	\|F_r'''\|_\infty\le12L+2M.
\end{equation}

For the hinge $h_k(x)=(x-k)_+$ in \eqref{eq:decomp}, set $A_k=\mathbb{E}(Z-k)_+$. Its Stein
solution is
\[
F_k(x)=
\begin{cases}
-A_k\,\dfrac{\Phi(x)}{\phi(x)},&x\le k,\\[1.2ex]
-1+(k+A_k)\dfrac{\overline{\Phi}(x)}{\phi(x)},&x>k.
\end{cases}
\]
For $|k|\le1$, both $A_k$ and $k+A_k$ lie in $[0,2]$. Moreover,
\[
\frac{\overline{\Phi}(x)}{\phi(x)}
=\int_0^\infty e^{-xu-u^2/2}\,du,
\qquad
\frac{\Phi(x)}{\phi(x)}
=\int_0^\infty e^{xu-u^2/2}\,du.
\]
On the corresponding sides of $k$, these representations yield
\[
|F_k'''(x)|
\le2\int_0^\infty u^3e^{u-u^2/2}\,du
\le2e\int_0^\infty u^3e^{-u^2/4}\,du
=16e.
\]
Combining this with \eqref{eq:sd-smooththird}, $|J_h|\le2L$, and
$F=J_hF_k+F_r$ proves the last bound in \eqref{eq:sd-steinregularity}.

\medskip
\noindent\emph{Step 2: Taylor expansion across the kink.}
Set $K=128(L+M)$. 
Since $F''$ has jump $J_h$ at $k$ and is $K$-Lipschitz on each side of $k$, we may write
\[
F''(t)=G(t)+J_h\mathbf 1_{\{t>k\}},
\]
where $G$ is globally $K$-Lipschitz. Hence, for $|t|\le |y|$,
\[
|F''(x+t)-F''(x)|
\le K|t|
+|J_h|\mathbf 1_{\{|x-k|\le |y|\}}.
\]
Integrating this estimate once gives
\begin{align}
	|F'(x+y)-F'(x)-yF''(x)|
	&\le\frac K2y^2
	+|J_h|\,|y|\mathbf{1}_{\{|x-k|\le|y|\}},\label{eq:sd-taylor1}
\end{align}
while integrating once more gives
\begin{align}
	\left|F(x+y)-F(x)-yF'(x)-\frac{y^2}{2}F''(x)\right|
	&\le\frac K6|y|^3
	+\frac{|J_h|}{2}y^2\mathbf{1}_{\{|x-k|\le|y|\}}.
	\label{eq:sd-taylor2}
\end{align}
The indicator accounts for the possible crossing of the kink $k$ by the segment joining $x$ and $x+y$.

Write
\[
Y^{(i)}=Y-Y_i,
\qquad
\sigma_i^2=\mathbb{E} Y_i^2,
\qquad
\mu_{3,i}=\mathbb{E} Y_i^3,
\qquad
p_i=\mathbb{P}(|Y^{(i)}-k|\le r_i).
\]
Define the Taylor remainders
\[
R_{1,i}
:=
F'(Y^{(i)}+Y_i)-F'(Y^{(i)})-Y_iF''(Y^{(i)})
\]
and
\[
R_{2,i}
:=
F(Y^{(i)}+Y_i)-F(Y^{(i)})
-Y_iF'(Y^{(i)})-\frac{Y_i^2}{2}F''(Y^{(i)}).
\]
Then, by independence of $Y_i$ and $Y^{(i)}$ and $\mathbb EY_i=0$,
\[
\mathbb EF'(Y)
=
\sum_i\sigma_i^2\mathbb EF'(Y^{(i)})
+\sum_i\sigma_i^2\mathbb ER_{1,i},
\]
where we used $\sum_i\sigma_i^2=1$. Moreover,
\[
\mathbb E[YF(Y)]
=
\sum_i\sigma_i^2\mathbb EF'(Y^{(i)})
+\frac12\sum_i\mu_{3,i}\mathbb EF''(Y^{(i)})
+\sum_i\mathbb E[Y_iR_{2,i}].
\]
Hence
\begin{equation}
\mathbb E[F'(Y)-YF(Y)]
=
-\frac12\sum_i\mu_{3,i}\mathbb EF''(Y^{(i)})+R, \label{eq:diff}
\end{equation}
where
\[
R:=
\sum_i\sigma_i^2\mathbb ER_{1,i}
-\sum_i\mathbb E[Y_iR_{2,i}].
\]

Indeed, \eqref{eq:sd-taylor1} gives
\[
|R_{1,i}|
\le \frac K2Y_i^2
+|J_h||Y_i|
\mathbf 1_{\{|Y^{(i)}-k|\le |Y_i|\}},
\]
and hence, using $|Y_i|\le r_i$, $\sigma_i^2\le r_i^2$, and independence,
\[
\sigma_i^2\mathbb E|R_{1,i}|
\le
\frac K2\sigma_i^4
+|J_h|r_i^3p_i.
\]
Similarly, \eqref{eq:sd-taylor2} yields
\[
|Y_iR_{2,i}|
\le
\frac K6|Y_i|^4
+\frac{|J_h|}{2}|Y_i|^3
\mathbf 1_{\{|Y^{(i)}-k|\le |Y_i|\}},
\]
so that
\[
\mathbb E|Y_iR_{2,i}|
\le
\frac K6\mathbb EY_i^4
+\frac{|J_h|}{2}r_i^3p_i.
\]
Summing these two estimates and recalling the definition of $R$ gives
\begin{equation}
|R|
\le
\frac K2\sum_i\sigma_i^4
+\frac K6\sum_i\mathbb EY_i^4
+\frac{3|J_h|}{2}\sum_i r_i^3p_i. \label{eq:sd-remainder}
\end{equation}


%
%

\medskip
\noindent\emph{Step 3: Controlling the kink contribution.}
Let
\[
s_i^2:=\operatorname{Var}(Y^{(i)})=1-\sigma_i^2\ge\frac34.
\]
By the Berry--Esseen bound,
\[
\sup_x\left|
\mathbb P(Y^{(i)}\le x)-\Phi(x/s_i)
\right|
\le
\frac{C_{\mathrm{BE}}\beta_3}{s_i^3}
\le 2C_{\mathrm{BE}}\beta_3.
\]
Therefore, applying the bound at the two endpoints \(k\pm r_i\),
\[
p_i
\le
\mathbb P\bigl(|s_iZ-k|\le r_i\bigr)
+4C_{\mathrm{BE}}\beta_3.
\]
Since \(s_i\ge\sqrt3/2\), the density of \(s_iZ\) is bounded by
$
1/(s_i\sqrt{2\pi}),
$
and hence
\[
\mathbb P\bigl(|s_iZ-k|\le r_i\bigr)
\le
\frac{2r_i}{s_i\sqrt{2\pi}}
\le r_i.
\]
Consequently,
\[
p_i\le r_i+4C_{\mathrm{BE}}\beta_3.
\]

Therefore,
\begin{align}
\sum_i r_i^3p_i
&\le\sum_i r_i^4
 +4C_{\mathrm{BE}}\beta_3\sum_i r_i^3\notag\\
&\le\sum_i r_i^4
 +4C_{\mathrm{BE}}\left(\sum_i r_i^3\right)^2\notag\\
&\le(1+4BC_{\mathrm{BE}})\sum_i r_i^4.
\label{eq:sd-kinkbound}
\end{align}
The last step is Cauchy--Schwarz:
\[
\left(\sum_i r_i^3\right)^2
\le\left(\sum_i r_i^2\right)\left(\sum_i r_i^4\right).
\]
Finally, substitute \eqref{eq:sd-kinkbound} into \eqref{eq:sd-remainder} and in turn into \eqref{eq:diff},
use $\|F''\|_\infty\le4L$ in \eqref{eq:sd-steinexpansion}, and apply
\eqref{eq:sd-stein}. This yields
\[
|\mathbb{E} h(Y)-\mathbb{E} h(Z)|
\le2L\sum_i|\mu_{3,i}|
 +\left(\frac{256}{3}+3+12BC_{\mathrm{BE}}\right)
  (L+M)\sum_i r_i^4,
\]
which implies \eqref{eq:sd-normal} with \eqref{eq:sd-CB}.
\end{IEEEproof}

\subsection{Exponential tilting and deficit estimates}

We use the standard exponential change of measure (Cram\'er/Esscher tilt);
see, e.g., \cite[Sec.~2.2]{DemboZeitouni2010}.
Let $S=\sum_i c_i X_i$ with $\sum_i c_i^2=1$.
By sign symmetry, assume $c_i\ge0$. Define
\[
K(t)=\log\mathbb{E} e^{tS}=\sum_i\log\cosh(tc_i)
\]
and
\begin{equation}\label{eq:sd-deficit}
\Delta=\Delta(t)=\frac{t^2}{2}-K(t)
=\sum_i\psi(tc_i),
\qquad
\psi(z)=\frac{z^2}{2}-\log\cosh z.
\end{equation}
For every $t>0$, one has $\Delta(t)>0$. Indeed, writing
\[
\psi(z):=\frac{z^2}{2}-\log\cosh z,
\]
we have
\[
\psi(0)=\psi'(0)=0,
\qquad
\psi''(z)=\tanh^2 z\ge0,
\]
with strict inequality for $z\neq0$. Hence
$
\psi(z)>0, z\neq0.
$
Since
\[
\Delta(t)=\sum_i\psi(tc_i)
\]
and $\sum_i c_i^2=1$, at least one $c_i$ is nonzero. Therefore, for every $t>0$, at least one term $\psi(tc_i)$ is strictly positive, and consequently
$
\Delta(t)>0.
$

For $t>0$, put
\[
g_t(y)=y_+e^{-ty},
\qquad
\mathcal{H}(t)=\mathbb{E} g_t(Z)=\int_0^\infty ye^{-ty}\phi(y)\,dy.
\]
A Gaussian change of variables gives
\begin{equation}\label{eq:sd-gaussiantilt}
\mathbb{E}(Z-t)_+=e^{-t^2/2}\mathcal{H}(t).
\end{equation}
For $t\ge1$, restriction to $0\le y\le1/t$ yields
\begin{equation}\label{eq:sd-calHlower}
\mathcal{H}(t)\ge\frac{c_*}{t^2},
\qquad
c_*:=\frac{e^{-1}\phi(1)}{2}>0.
\end{equation}

Introduce the tilted measure
\[
\frac{d\mathbb{P}_t}{d\mathbb{P}}=e^{tS-K(t)}.
\]
Then
\begin{equation}\label{eq:sd-tiltidentity}
\mathbb{E}(S-t)_+
=e^{-t^2/2-\Delta}\mathbb{E}_tg_t(S-t).
\end{equation}
Under $\mathbb{P}_t$, the signs remain independent, and
\[
\mathbb{E}_t X_i=\tanh(tc_i) 
\qquad
\operatorname{Var}_t(X_i)
=1-\tanh^2(tc_i)
=\operatorname{sech}^2(tc_i).
\]
Write
\begin{equation}\label{eq:sd-mvd}
m=K'(t),\qquad
v^2=K''(t)=\sum_i c_i^2\operatorname{sech}^2(tc_i),\qquad
d=t-m,\qquad
q_4=\sum_i c_i^4.
\end{equation}
In particular, $0<v\le1$ and $d\ge0$.

\begin{lem}[Deficit controls]
\label{lem:tilt-deficit}
For every $t>0$,
\begin{align}
1-v^2&\le\frac{48\Delta}{t^2},
& 0\le d&\le\frac{48\Delta}{t},\label{eq:sd-vdcontrol}\\
q_4&\le\frac{1152\Delta(1+\Delta)}{t^4},
&\max_i c_i^2&\le\frac{1+48\Delta}{t^2}.
\label{eq:sd-qcontrol}
\end{align}
\end{lem}

\begin{IEEEproof}
	Set
$
	z_i:=tc_i,\, i\in[n].
$
	Recall that
	\[
	\Delta=\sum_i\psi(z_i),\qquad
	\psi(z):=\frac{z^2}{2}-\log\cosh z,
	\]
	and that
	\[
	v^2=\sum_i c_i^2\operatorname{sech}^2 z_i,\qquad
	d=t-\sum_i c_i\tanh z_i,\qquad
	\sum_i c_i^2=1.
	\]
	
	We first establish a lower bound on $\psi$. For $z\ge0$,
	\[
	\tanh z=\frac{e^{2z}-1}{e^{2z}+1}
	\ge \frac{2z}{2+2z}
	=\frac{z}{1+z},
	\]
	where we used $e^{2z}\ge1+2z$. Since
	\[
	\psi'(z)=z-\tanh z,
	\qquad
	\psi''(z)=\tanh^2z,
	\qquad
	\psi(0)=\psi'(0)=0,
	\]
	we have
	\[
	\psi(z)
	=\int_0^z(z-u)\tanh^2u\,du.
	\]
	Using $\tanh u\ge u/(1+u)$ and $u\le z$,
	\begin{align}
		\psi(z)
		&\ge
		\frac{1}{(1+z)^2}
		\int_0^z(z-u)u^2\,du 
		=
		\frac{z^4}{12(1+z)^2}
		\ge
		\frac{z^4}{24(1+z^2)},
		\label{eq:sd-psilower}
	\end{align}
	where the last inequality follows from
	$(1+z)^2\le2(1+z^2)$.
	
	We next note that
	\begin{equation}
		z^2\tanh^2z
		\le \frac{2z^4}{1+z^2}
		\le48\psi(z).
		\label{eq:sd-tanhcontrol}
	\end{equation}
	Indeed, the first inequality follows from $\tanh z\le z$ when
	$0\le z\le1$, and from $\tanh z\le1$ when $z\ge1$; the second follows
	from \eqref{eq:sd-psilower}. Moreover,
	\[
	z-\tanh z
	=\int_0^z\tanh^2u\,du,
	\]
	and hence, since $\tanh$ is increasing,
	\[
	z(z-\tanh z)
	=
	z\int_0^z\tanh^2u\,du
	\le z^2\tanh^2z.
	\]
	
	We can now control $v$ and $d$. Since
	$1-\operatorname{sech}^2z=\tanh^2z$,
	\begin{align*}
		1-v^2
		&=\sum_i c_i^2\tanh^2 z_i =\frac1{t^2}\sum_i z_i^2\tanh^2z_i
		\le\frac{48}{t^2}\sum_i\psi(z_i)
		=\frac{48\Delta}{t^2}.
	\end{align*}
	Likewise, using $t=t\sum_i c_i^2=\sum_i c_i z_i$,
	\begin{align*}
		d
		&=t-\sum_i c_i\tanh z_i\\
		&=\sum_i c_i(z_i-\tanh z_i)\\
		&=\frac1t\sum_i z_i(z_i-\tanh z_i)
		\le\frac1t\sum_i z_i^2\tanh^2z_i
		\le\frac{48\Delta}{t}.
	\end{align*}
	
	Furthermore, since
	$
	d:=t-m,\,
	m=\sum_i c_i\tanh(tc_i),
	$
	and $c_i\ge0$, the inequality $\tanh z\le z$ for $z\ge0$ gives
	\[
	m
	\le
	\sum_i c_i(tc_i)
	=
	t\sum_i c_i^2
	=
	t.
	\]
	Therefore
	$
	d=t-m\ge0.
	$
	This proves \eqref{eq:sd-vdcontrol}.
	
	It remains to control the largest coefficient and $q_4$. Let
	\[
	q:=\max_i z_i^2=\max_i (tc_i)^2.
	\]
	Choose $j$ with $z_j^2=q$. By \eqref{eq:sd-psilower},
	\[
	\Delta\ge\psi(z_j)
	\ge\frac{q^2}{24(1+q)},
	\]
	so
	\[
	\frac{q^2}{1+q}\le24\Delta.
	\]
	If $q\le1$, then trivially $q\le1+48\Delta$. If $q>1$, then
	\[
	\frac{q^2}{1+q}\ge\frac q2,
	\]
	and therefore $q\le48\Delta$. Thus, in all cases,
$
	q\le1+48\Delta.
$
	Consequently,
	\[
	\max_i c_i^2=\frac{q}{t^2}
	\le\frac{1+48\Delta}{t^2}.
	\]
	
	Finally,
	\begin{align*}
		t^4q_4
		&=\sum_i z_i^4 =\sum_i(1+z_i^2)\frac{z_i^4}{1+z_i^2} \le(1+q)\sum_i\frac{z_i^4}{1+z_i^2}\\
		&\le24(1+q)\sum_i\psi(z_i) \le24(2+48\Delta)\Delta \le1152\Delta(1+\Delta),
	\end{align*}
	where we used $1+q\le2+48\Delta$ and
	$24(2+48\Delta)\le1152(1+\Delta)$. Dividing by $t^4$ gives
	\[
	q_4\le\frac{1152\Delta(1+\Delta)}{t^4},
	\]
	which, together with the preceding bound on $\max_i c_i^2$, proves
	\eqref{eq:sd-qcontrol}.
\end{IEEEproof}

%

\subsection{Proof of eventual Gaussian domination}

\begin{IEEEproof}[Proof of Theorem~\ref{thm:stoploss}]
Let
\[
C_8:=100+96C_{\mathrm{BE}}.
\]
Indeed, in the application below (with $r_i$ chosen in \eqref{eq:ri}),
\[
\sum_i r_i^2=\frac{4}{v^2}\le 8
\]
by \eqref{eq:sd-admissible}. Hence Lemma~\ref{lem:normal-kink}, for which
\[
C_B=100+12BC_{\mathrm{BE}},
\]
is applied with $B=8$, giving
\[
C_8=100+12\cdot 8\,C_{\mathrm{BE}}
=100+96C_{\mathrm{BE}}.
\]

Define
\begin{equation}\label{eq:sd-constants}
A=200C_8,
\qquad
C_*=\frac{1152A}{c_*}.
\end{equation}
One sufficient, deliberately nonoptimized choice of threshold is
\begin{equation}\label{eq:sd-Tchoice}
T_0=\max\left\{
10^6,\ 144(1+C_*)^2,\ \frac{2}{ec_*}
\right\}.
\end{equation}
Fix $t\ge T_0$ and abbreviate $\Delta=\Delta(t)>0$.
Equations \eqref{eq:sd-gaussiantilt} and \eqref{eq:sd-tiltidentity} imply
\begin{equation}\label{eq:sd-ratio}
\frac{\mathbb{E}(S-t)_+}{\mathbb{E}(Z-t)_+}
=e^{-\Delta}\frac{\mathbb{E}_tg_t(S-t)}{\mathcal{H}(t)}.
\end{equation}

\medskip
\noindent\emph{Case 1: $\Delta\ge2\log t$.}
Since $0\le g_t(y)\le1/(et)$ (from the definition of $g_t$), \eqref{eq:sd-calHlower} gives
\[
\frac{\mathbb{E}(S-t)_+}{\mathbb{E}(Z-t)_+}
\le e^{-\Delta}\frac{t}{ec_*}
\le\frac{1}{ec_*t}
\le\frac12<1.
\]

\medskip
\noindent\emph{Case 2: $0<\Delta\le2\log t$.}
The deficit controls imply
\begin{equation}\label{eq:sd-admissible}
v^2\ge\frac12,
\qquad
\frac dv\le1,
\qquad
\max_i\frac{2c_i}{v}\le\frac12.
\end{equation}
These conclusions already hold for $t\ge10^6$:
using $\log t\le\sqrt t$ for $t\ge1$, we have
\[
1-v^2\le\frac{96\log t}{t^2},
\qquad
\frac dv\le\frac{96\sqrt2\log t}{t},
\qquad
\max_i\frac{4c_i^2}{v^2}
\le\frac{8(1+96\log t)}{t^2}.
\]
The first of these bounds gives $v^2\ge1/2$ before the other two
are applied.

Under $\mathbb{P}_t$, define
\begin{equation}\label{eq:ri}
Y_i=\frac{c_i( X_i-\tanh(tc_i))}{v},
\qquad
r_i=\frac{2c_i}{v}.
\end{equation}
Then the $Y_i$ are independent and centered, and
\[
\sum_i\mathbb{E}_tY_i^2=1,
\qquad
|Y_i|\le r_i\le\frac12,
\qquad
\sum_i r_i^2=\frac4{v^2}\le8.
\]
For a sign variable of mean $\theta$,
\[
\mathbb{E}(X-\theta)^3=-2\theta(1-\theta^2).
\]
Consequently, since $\tanh(tc_i)\le tc_i$,
\begin{equation}\label{eq:sd-tiltedmoments}
\sum_i|\mathbb{E}_tY_i^3|
\le\frac{2t}{v^3}q_4\le6tq_4,
\qquad
\sum_i r_i^4=\frac{16}{v^4}q_4\le64q_4.
\end{equation}

Apply Lemma~\ref{lem:normal-kink} to
\[
h(x)=g_t(vx-d).
\]
Its kink is at $d/v\in[0,1]$, and
\[
\operatorname{Lip}(h)\le1,
\qquad
|h''(x)|\le2t\quad(x\ne d/v).
\]
Since by definition, $S-t=v\sum_iY_i-d$, \eqref{eq:sd-tiltedmoments} yields
\begin{align}
\left|\mathbb{E}_tg_t(S-t)-\mathbb{E} g_t(vZ-d)\right|
&\le C_8\bigl[6t+64(1+2t)\bigr]q_4\notag\\
&\le Atq_4.
\label{eq:sd-tiltedapprox}
\end{align}
Because $d\ge0$ and $v\le1$,
\begin{align}
\mathbb{E} g_t(vZ-d)
&=\frac1v\int_0^\infty
 ye^{-ty}\phi\!\left(\frac{y+d}{v}\right)\,dy\notag\\
&\le\frac1v\int_0^\infty ye^{-ty}\phi(y)\,dy
=\frac{\mathcal{H}(t)}v.
\label{eq:sd-gaussiancomparison}
\end{align}
Combining \eqref{eq:sd-ratio}, \eqref{eq:sd-tiltedapprox}, and
\eqref{eq:sd-gaussiancomparison},
\begin{equation}\label{eq:sd-ratio2}
\frac{\mathbb{E}(S-t)_+}{\mathbb{E}(Z-t)_+}
\le e^{-\Delta}\left[\frac1v+\frac{Atq_4}{\mathcal{H}(t)}\right].
\end{equation}
For $0\le u\le1/2$, $(1-u)^{-1/2}\le1+2u$. Hence
\eqref{eq:sd-vdcontrol} gives
\[
\frac1v\le1+\frac{96\Delta}{t^2}.
\]
Using \eqref{eq:sd-calHlower} and \eqref{eq:sd-qcontrol} in
\eqref{eq:sd-ratio2}, we obtain
\begin{align}
\frac{\mathbb{E}(S-t)_+}{\mathbb{E}(Z-t)_+}
&\le e^{-\Delta}\left[
1+\frac{96\Delta}{t^2}
+\frac{C_*\Delta(1+\Delta)}{t}
\right]\notag\\
&\le e^{-\Delta}\left[
1+\Delta\left(
\frac{96}{t^2}+\frac{C_*(1+2\log t)}{t}
\right)
\right].
\label{eq:sd-finalratio}
\end{align}
The choice \eqref{eq:sd-Tchoice}, together with
$1+2\log t\le3\sqrt t$ for $t\ge1$, ensures
\[
\frac{96}{t^2}+\frac{C_*(1+2\log t)}{t}
\le\frac{96}{t^2}+\frac{3C_*}{\sqrt t}
\le\frac12.
\]
Therefore
\[
\frac{\mathbb{E}(S-t)_+}{\mathbb{E}(Z-t)_+}
\le e^{-\Delta}\left(1+\frac\Delta2\right)<1,
\]
since $\Delta>0$. This proves the desired strict inequality.
\end{IEEEproof}

\subsection{Application to level-$1$ Fourier weight}

\begin{IEEEproof}[Proof of Theorem~\ref{thm:small-density}]
We first prove the upper bound at dyadic densities. Let
$f:\{\pm1\}^n\to\{0,1\}$ have mean $a\in(0,a_0]$, and put
\[
b_i=\hat f_{i},
\qquad
b=\left(\sum_i b_i^2\right)^{1/2}.
\]
If $b=0$, then $\mathbf{W}_{1}[f]<J(a)$ is immediate. Otherwise define
\[
S=\sum_i\frac{b_i}{b} X_i.
\]
Then $\operatorname{Var}(S)=1$ and $\mathbb{E}[fS]=b$. For every $t\in\mathbb{R}$,
\begin{equation}\label{eq:sd-fourierdual}
b=ta+\mathbb{E}[f(S-t)]\le ta+\mathbb{E}(S-t)_+.
\end{equation}
Take $t=\overline{\Phi}^{-1}(a)\ge T_0$. By Theorem~\ref{thm:stoploss},
\[
b<ta+\mathbb{E}(Z-t)_+
=ta+\phi(t)-t\overline{\Phi}(t)
=\phi(t).
\]
Thus $\mathbf{W}_{1}[f]<J(a)$, and taking the supremum over dimensions gives
$W(a)\le J(a)$.

For the reverse bound, fix any dyadic $a\in(0,1)$. For sufficiently
large $n$, let $A_n$ consist of the $a2^n$ vertices with largest values
of
\[
S_n=\frac{X_1+\cdots+X_n}{\sqrt n},
\]
using an arbitrary subset of the last layer when needed, and set
$f_n=\mathbf{1}_{A_n}$. Then Cauchy--Schwarz yields
\begin{equation}\label{eq:sd-lowerfourier}
\mathbf{W}_{1}[f_n]
\ge\left(\frac1{\sqrt n}\sum_i\widehat{f_n}_i\right)^2
=\bigl(\mathbb{E}[f_nS_n]\bigr)^2.
\end{equation}
The central limit theorem and continuity of the Gaussian distribution
imply convergence of the upper-$a$ quantiles of $S_n$ to
$t=\overline{\Phi}^{-1}(a)$. Also, $\sup_n\mathbb{E} S_n^2=1$ implies uniform
integrability. These are standard consequences of weak convergence and
uniform integrability; see, e.g., \cite[Chs.~2 and~5]{Billingsley1999}. Therefore
\[
\mathbb{E}[f_nS_n]
\longrightarrow\mathbb{E}[Z\mathbf{1}_{\{Z>t\}}]
=\phi(t).
\]
The partially filled threshold layer causes no difficulty: its probability
tends to zero, for example by \eqref{eq:sd-BE}. Equation
\eqref{eq:sd-lowerfourier} now gives $W(a)\ge J(a)$.

For completeness, we also justify the continuous extension of $W$ used throughout the paper. If $a,b\in\mathcal D$,
embed a Boolean function of mean $a$ into a sufficiently large cube
and add or remove exactly a fraction $|a-b|$ of the vertices to obtain
a Boolean function $g$ of mean $b$. Parseval's identity gives
\[
\left\|(\hat f_{i})_{i\in [n]}-(\hat g_i)_{i\in [n]}\right\|_2
\le\|f-g\|_2=\sqrt{|a-b|}.
\]
Taking suprema and interchanging $a,b$ yields
\begin{equation}\label{eq:sd-continuity}
|\sqrt{W(a)}-\sqrt{W(b)}|\le\sqrt{|a-b|}.
\end{equation}
Thus $W$ has a unique continuous extension to $[0,1]$. Since $J$ is
continuous, equality extends from dyadic densities to every real
$a\in(0,a_0]$.
\end{IEEEproof}

%

The exact theorem above is strongest in the extreme small-density
regime. We next develop explicit upper bounds that remain informative
throughout the full interval $a\in(0,1/2]$.

\section{Improved Bounds on Level-$1$ Weight}\label{sec:improved-bounds}

\subsection{Improved Bound for $a<1/4$}

We now improve the level-$1$ weight bound in Theorem~\ref{thm:level-1-bound}
by using an inductive method.   
\begin{prop}
\label{prop:bound}Let $\varphi$ be an upper bound on $W$, i.e.,
$W(a)\le\varphi(a),\forall a\in[0,1].$ Suppose $\chi:[0,1]\to[0,\infty)$
is a function such that 
\begin{equation}
\chi(0)=\chi(1)=0,\label{eq:-4}
\end{equation}
and
\begin{equation}
\Gamma(a_{0},a_{1})\le0,\forall a_{0},a_{1}\in[0,1],\label{eq:-6}
\end{equation}
where 
\begin{equation}
\Gamma(a_{0},a_{1}):=\min\left\{ \frac{1}{4}\left(\sqrt{\chi(a_{0})}+\sqrt{\chi(a_{1})}\right)^{2}+\frac{1}{4}(a_{1}-a_{0})^{2},\varphi\left(\frac{a_{0}+a_{1}}{2}\right)\right\} -\chi\left(\frac{a_{0}+a_{1}}{2}\right).\label{eq:-6-1}
\end{equation}
 Then, it holds that 
\[
W(a)\le\chi(a),\forall a\in[0,1].
\]
\end{prop}
\begin{IEEEproof}
For $n=0$, $W^{(0)}(0)=\chi(0)=0,W^{(0)}(1)=\chi(1)=0$. Assume
for $n=k$, $W^{(k)}(a)\le\chi(a)$. For $n=k+1$, let $f$ be a Boolean
function attaining $W^{(k+1)}(a)$. We can write $f(x_{1},x_{2},...,x_{k+1})=\frac{1+x_{k+1}}{2}g(x_{1},\ldots,x_{k})+\frac{1-x_{k+1}}{2}h(x_{1},\ldots,x_{k})$,
where $g(x_{1},\ldots,x_{k})=f(x_{1},...,x_{k},1)$ and $h(x_{1},\ldots,x_{k})=f(x_{1},...,x_{k},-1)$.
Then, $\hat{f}_{i}=\frac{1}{2}\hat{g}_{i}+\frac{1}{2}\hat{h}_{i},\forall i\in[k]$,
and $\hat{f}_{k+1}=\frac{1}{2}\mathbb{E}g-\frac{1}{2}\mathbb{E}h$.
Moreover, $a=\frac{1}{2}\mathbb{E}g+\frac{1}{2}\mathbb{E}h$. So,
$\mathbb{E}g=a+\hat{f}_{k+1}$ and $\mathbb{E}h=a-\hat{f}_{k+1}$.
Note that $g$ and $h$ are Boolean functions on the $k$-dimensional
space. Hence, their Fourier coefficients $\hat{g}_{i},\hat{h}_{i},\forall i\in[k]$
are also defined on the $k$-dimensional space. By the Minkowski inequality,
\begin{align}
W^{(k+1)}(a) & =\sum_{i=1}^{k+1}\hat{f}_{i}^{2}=\frac{1}{4}\sum_{i=1}^{k}(\hat{g}_{i}+\hat{h}_{i})^{2}+\hat{f}_{k+1}^{2}\nonumber \\
 & \leq\frac{1}{4}\Big(\sqrt{\sum_{i=1}^{k}\hat{g}_{i}^{2}}+\sqrt{\sum_{i=1}^{k}\hat{h}_{i}^{2}}\Big)^{2}+\hat{f}_{k+1}^{2}\nonumber \\
 & =\frac{1}{4}\left(\sqrt{\mathbf{W}_{1}[g]}+\sqrt{\mathbf{W}_{1}[h]}\right)^{2}+\hat{f}_{k+1}^{2}\nonumber \\
 & \le\frac{1}{4}\left(\sqrt{\chi(a+\hat{f}_{k+1})}+\sqrt{\chi(a-\hat{f}_{k+1})}\right)^{2}+\hat{f}_{k+1}^{2}.\label{eq:}
\end{align}
Moreover, $W^{(k+1)}(a)\le\varphi\left(a\right)$. Hence, 
\begin{align*}
W^{(k+1)}(a) & \le\min\left\{ \frac{1}{4}\left(\sqrt{\chi(a+\hat{f}_{k+1})}+\sqrt{\chi(a-\hat{f}_{k+1})}\right)^{2}+\hat{f}_{k+1}^{2},\varphi\left(a\right)\right\} \le\chi\left(a\right),
\end{align*}
where the last inequality follows by the assumption in \eqref{eq:-6}. 
\end{IEEEproof}

In \eqref{eq:-6}, we require $\chi$ to satisfy the inequality for
all $a_{0},a_{1}\in[0,1]$. We now focus on a specific form of $\chi$
and for this kind of $\chi$, we relax the requirement in \eqref{eq:-6}.
Let $T\in(0,1/4]$. We now consider 
\begin{equation}
\chi(a)=\begin{cases}
I_{w}^{2}(a), & a\in[0,T]\\
\varphi\left(a\right), & a\in(T,1/2]
\end{cases},\label{eq:J}
\end{equation}
where 
\[
I_{w}(a):=wI\left(\frac{a}{w}\right)
\]
with $w$ denoting the unique solution to 
\begin{equation}
I_{w}(T)=\sqrt{\varphi(T)}.\label{eq:w}
\end{equation}
So, $\chi$ is continuous (when $\varphi$ is continuous). Moreover,
it is easy to see that given any $v>0$, it always holds that $I_{v}(a)I_{v}''(a)=-1$. 

\begin{prop}
\label{prop:bound2} If  
\begin{equation}
\Gamma(a_{0},a_{1})\le0,\;\forall a_{0},a_{1}\in[0,T+\sqrt{\frac{T}{6}}],\frac{a_{0}+a_{1}}{2}\le T,\label{eq:-9}
\end{equation}
then for $a\in[0,1/2]$, 
\[
W(a)\le\chi(a).
\]
\end{prop}
\begin{IEEEproof}
For $n\le2$, $a$ only takes values in $\{0,1/4,1/2\}$. For these
values, $W(0)=\chi(0)=0,$$W(1/4)\le\varphi\left(1/4\right)=\chi(1/4),$
and $W(1/2)\le\varphi\left(1/2\right)=\chi(1/2)$. So, for $n\le2$,
it holds that $W^{(n)}(a)\le\chi(a).$ 

We next consider the case $n\ge3$. Since $\sum_{i=1}^{n}\hat{f}_{i}^{2}\le W(a)$,
we have that there is some $i$ such that $|\hat{f}_{i}|\le\sqrt{\frac{W(a)}{n}}$.
By symmetry, we assume $|\hat{f}_{n}|\le\sqrt{\frac{W(a)}{n}}$. So,
both $a+\hat{f}_{n}$ and $a-\hat{f}_{n}$ are upper bounded by $a+\sqrt{\frac{W(a)}{n}}$.
Since $W(a)\le a/2$ for $a\le1/2$, we have that for $a\le T$ and
$n\ge3$, $a+\sqrt{\frac{W(a)}{n}}\le T+\sqrt{\frac{T}{6}}$.  

For $n\ge3$ and $a\in[T,\frac{1}{2}]$, it holds that $W^{(n)}(a)\le\varphi\left(a\right)=\chi\left(a\right)$.
For $n\ge3$ and $a\in[0,T]$, by the induction in the proof of Proposition
\ref{prop:bound} but with the assumption in \eqref{eq:-6} replaced
by the one in \eqref{eq:-9}, we can verify that $W^{(n)}(a)\le\chi(a)$. 
\end{IEEEproof}

We choose for $0\leq t\leq1$, 
\[
\varphi(t)=\tilde{t}/2,
\]
where $\tilde{t}:=\min\{t,1-t\}$.  For such $\varphi$, we choose
$T=0.21$, which implies that the unique solution to \eqref{eq:w}
satisfies $w\approx1.36971$. We now state our main result. 
\begin{thm}
\label{thm:bound}The function $\chi$ given in \eqref{eq:J} with
$\varphi(t)=t/2$ and $T=0.21$ satisfies the conditions in Proposition
\ref{prop:bound2}. As a consequence, it holds that for $a\in[0,1/2]$,
\[
W(a)\le\chi(a).
\]
\end{thm}
\begin{IEEEproof}
Observe that $T+\sqrt{\frac{T}{6}}\le0.4$. So, we only need prove
that 
\begin{equation}
\Gamma(a_{0},a_{1})\le0,\forall a_{0},a_{1}\in[0,0.4],\frac{a_{0}+a_{1}}{2}\le T.\label{eq:-10}
\end{equation}

We first prove this inequality for $a_{0},a_{1}\in[0,T]$, i.e., 
\[
\Gamma(a_{0},a_{1})=\frac{1}{4}\left(I_{w}(a_{0})+I_{w}(a_{1})\right)^{2}+\frac{1}{4}(a_{1}-a_{0})^{2}-I_w^{2}\left(\frac{a_{0}+a_{1}}{2}\right)\le0.
\]
We compute the derivative 
\begin{align*}
\Gamma_{2}'(a_{0},a_{1}):=\partial_{a_{1}}\Gamma(a_{0},a_{1}) & =\frac{1}{2}\left(I_{w}(a_{0})+I_w(a_{1})\right)I_w'(a_{1})+\frac{1}{2}(a_{1}-a_{0})\\
 & \qquad-I_w\left(\frac{a_{0}+a_{1}}{2}\right)I_w'\left(\frac{a_{0}+a_{1}}{2}\right),
\end{align*}
and the second derivative 
\begin{align*}
\Gamma_{1,2}''(a_{0},a_{1}):=\partial_{a_{0}}\partial_{a_{1}}\Gamma(a_{0},a_{1}) & =\frac{1}{2}I_w'(a_{0})I_w'(a_{1})-\frac{1}{2}-\frac{1}{2}\left(I_w'\left(\frac{a_{0}+a_{1}}{2}\right)\right)^{2}\\
 & \qquad-\frac{1}{2}I_w\left(\frac{a_{0}+a_{1}}{2}\right)I_w''\left(\frac{a_{0}+a_{1}}{2}\right).
\end{align*}
Since $I_w(a)I_w''(a)=-1$, denoting $g(a):=I_{w}'(a)=-\Phi^{-1}\left(\frac{a}{w}\right)$,
we have that 
\begin{align*}
\Gamma_{1,2}'' & (a_{0},a_{1})=\frac{1}{2}g(a_{0})g(a_{1})-\frac{1}{2}g^{2}\left(\frac{a_{0}+a_{1}}{2}\right).
\end{align*}

We claim that $g$ is log-convex on $[0,T]$, i.e., $g''g\ge(g')^{2}$,
which implies $\Gamma_{1,2}''(a_{0},a_{1})\ge0$ for $a_{0},a_{1}\in[0,T]$.
We now prove this claim. Observe that 
\begin{align}
g(a) & =-\Phi^{-1}\left(\frac{a}{w}\right)\nonumber \\
g'(a) & =\frac{-1}{w\phi(g(a))}\label{eq:-16}\\
g''(a) & =\frac{\phi'(g(a))g'(a)}{w\phi(g(a))^{2}}=-\frac{g(a)\phi(g(a))g'(a)}{w\phi(g(a))^{2}}\nonumber \\
 & =-\frac{g(a)g'(a)}{w\phi(g(a))}=g(a)\left(g'(a)\right)^{2}.\label{eq:-5}
\end{align}
So, $g''g-(g')^{2}=g^{2}\left(g'\right)^{2}-(g')^{2}=(g^{2}-1)(g')^{2}$.
We only need prove $g\ge1$ on $[0,T]$. Since $g$ is decreasing
on $[0,1]$, it suffices to prove $g(T)\ge1$. Noting $T=0.21$, we
can verify that $g(T)=1.02231\ldots\ge1$. This completes the proof of
the claim above. 

By the claim above, $\Gamma_{1,2}''(a_{0},a_{1})\ge0$ for $a_{0},a_{1}\in[0,T]$.
We assume $a_{0}\le a_{1}$. So, $\Gamma_{2}'(a_{0},a_{1})\le\Gamma_{2}'(a_{1},a_{1})$
for $a_{0},a_{1}\in[0,T]$. Furthermore,
\[
\Gamma_{2}'(a_{1},a_{1})=I_{w}\left(a_{1}\right)I_{w}'(a_{1})-I_{w}\left(a_{1}\right)I_{w}'\left(a_{1}\right)=0.
\]
So, 
\begin{equation}
\Gamma_{2}'(a_{0},a_{1})\le0\label{eq:-11}
\end{equation}
for $0\le a_{0}\le a_{1}\le T$, which implies that given $a_{0}$,
$\Gamma(a_{0},a_{1})$ is decreasing in $a_{1}\in[a_{0},T]$. So,
$\Gamma(a_{0},a_{1})\le\Gamma(a_{0},a_{0})=0$ for $a_{0},a_{1}\in[0,T]$. 

We next prove the inequality in \eqref{eq:-10} for $a_{1}>T,a_{0}+a_{1}\le2T$,
i.e., 
\[
\Gamma(a_{0},a_{1})=\frac{1}{4}\left(I_{w}(a_{0})+\sqrt{\frac{a_{1}}{2}}\right)^{2}+\frac{1}{4}(a_{1}-a_{0})^{2}-I_{w}^{2}\left(\frac{a_{0}+a_{1}}{2}\right)\le0.
\]
We compute derivatives of $\Gamma$: 
\[
\Gamma_{2}'(a_{0},a_{1})=\frac{1}{2}\left(I_{w}(a_{0})+\sqrt{\frac{a_{1}}{2}}\right)\frac{1}{2\sqrt{2a_{1}}}+\frac{1}{2}(a_{1}-a_{0})-I_{w}\left(\frac{a_{0}+a_{1}}{2}\right)I_{w}'\left(\frac{a_{0}+a_{1}}{2}\right),
\]
\begin{align*}
\Gamma_{2,2}''(a_{0},a_{1}) & =\frac{1}{2}I_{w}(a_{0})\frac{-1}{4\sqrt{2a_{1}^{3}}}+\frac{1}{2}-\frac{1}{2}\left(I_{w}'\left(\frac{a_{0}+a_{1}}{2}\right)\right)^{2}-\frac{1}{2}I_{w}\left(\frac{a_{0}+a_{1}}{2}\right)I_{w}''\left(\frac{a_{0}+a_{1}}{2}\right)\\
 & =1-\frac{1}{2}\frac{I_{w}(a_{0})}{4\sqrt{2a_{1}^{3}}}-\frac{1}{2}\left(I_{w}'\left(\frac{a_{0}+a_{1}}{2}\right)\right)^{2},
\end{align*}
and 
\begin{align*}
\Gamma_{2,2,2}'''(a_{0},a_{1}) & =\frac{1}{2}\frac{3I_{w}(a_{0})}{8\sqrt{2a_{1}^{5}}}-\frac{1}{2}I_{w}'\left(\frac{a_{0}+a_{1}}{2}\right)I_{w}''\left(\frac{a_{0}+a_{1}}{2}\right)\\
 & =\frac{3I_{w}(a_{0})}{8\sqrt{2a_{1}^{5}}}+\frac{I_{w}'\left(\frac{a_{0}+a_{1}}{2}\right)}{I_{w}\left(\frac{a_{0}+a_{1}}{2}\right)}\ge0,
\end{align*}
where the last inequality follows since $I_{w}'\left(\frac{a_{0}+a_{1}}{2}\right)=g(\frac{a_{0}+a_{1}}{2})\ge g(T)=1.02231\ldots\ge1$
due to the fact that $\frac{a_{0}+a_{1}}{2}\ge T$ and $g$ is decreasing.
 So, $\Gamma_{2,2}''(a_{0},a_{1})$ is increasing in $a_{1}$. Moreover,
observe that
\begin{align*}
\Gamma_{2,2}''(a_{0},2T-a_{0}) & =1-\frac{1}{2}\frac{I_{w}(a_{0})}{4\sqrt{2a_{1}^{3}}}-\frac{1}{2}\left(I_{w}'\left(\frac{a_{0}+a_{1}}{2}\right)\right)^{2}\\
 & =1-\frac{1}{2}\frac{I_{w}(a_{0})}{4\sqrt{2(2T-a_{0})^{3}}}-\frac{1}{2}\left(I_{w}'\left(T\right)\right)^{2}\\
 & \ge1-\frac{1}{2}\frac{I_{w}(T)}{4\sqrt{2(2T-T)^{3}}}-\frac{1}{2}\left(I_{w}'\left(T\right)\right)^{2}\\
 & =1-\frac{1}{2}\frac{I_{w}(T)}{4\sqrt{2T^{3}}}-\frac{1}{2}g(T)^{2}\\
 & =0.179822\ldots\ge0.
\end{align*}
So, $\Gamma_{2,2}''(a_{0},a_{1})$ is either positive or first-negative-then-positive
in $a_{1}\in[T,2T-a_{0}]$. So, $\Gamma_{2}'(a_{0},a_{1})$ is either
increasing or first-decreasing-then-increasing in $a_{1}\in[T,2T-a_{0}]$.
We now require the following lemma, the proof of which is given in
Appendix~\ref{sec:Proof-of-Lemma}. 
\begin{lem}
\label{lem:Gamma}It holds that $\lim_{a_{1}\downarrow T}\Gamma_{2}'(a_{0},a_{1})\le0$.
Moreover, $\Gamma_{2}'(a_{0},2T-a_{0})\le0$ for $0.02\le a_{0}\le T$,
and $\Gamma_{2}'(a_{0},0.4)\le0$ for $a_{0}\le0.02$.
\end{lem}
Combining this lemma with the fact that $\Gamma_{2}'(a_{0},a_{1})$
is either increasing or first-decreasing-then-increasing in $a_{1}\in[T,\min\{2T-a_{0},0.4\}]$,
we obtain that $\Gamma_{2}'(a_{0},a_{1})\le0$ for $a_{1}\in[T,\min\{2T-a_{0},0.4\}]$.
So, $\Gamma(a_{0},a_{1})$ is decreasing in $a_{1}\in[T,\min\{2T-a_{0},0.4\}]$.
So, $\Gamma(a_{0},a_{1})\le\Gamma(a_{0},T)\le0,$ where the last inequality
was proven in the first case. 
\end{IEEEproof}

\subsection{Improved Bound for $a>1/4$}

We now improve Chang's lemma for $a>1/4$. 
\begin{thm}
\label{thm:bound2}For $a\in[0,1/2]$,
\begin{align}
W(a)\le\tilde{\chi}(a) & :=\max_{\beta\in[0,a]}\min\Big\{\frac{1}{4}\left(\sqrt{\chi(a+\beta)}+\sqrt{\chi(a-\beta)}\right)^{2}+\beta^{2},\nonumber \\
 & \qquad\frac{1}{4}\left(\sqrt{4\left(\frac{1}{2}-\frac{1}{\sqrt{2\pi}}\right)\beta+\frac{1}{2\pi}}+\frac{1}{\sqrt{2\pi}}\right)^{2}\Big\}.\label{eq:-19}
\end{align}
\end{thm}
\begin{IEEEproof}
Let $f$ be a Boolean function attaining $W^{(n)}(a)$. Let $\beta=\max_{i\in[n]}|\hat{f}_{i}|$.
Then, from \eqref{eq:}, it holds that 
\begin{align}
W^{(n)}(a) & \le\frac{1}{4}\left(\sqrt{W^{(n-1)}(a+\beta)}+\sqrt{W^{(n-1)}(a-\beta)}\right)^{2}+\beta^{2}.\label{eq:-1}
\end{align}
Applying $W^{(n-1)}(t)\le\chi(t)$, we obtain 
\begin{align}
W(a) & \le\frac{1}{4}\left(\sqrt{\chi(a+\beta)}+\sqrt{\chi(a-\beta)}\right)^{2}+\beta^{2}.\label{eq:b1}
\end{align}

On the other hand, we claim that 
\begin{align}
W(a) & \le\frac{1}{4}\left(\sqrt{4\left(\frac{1}{2}-\frac{1}{\sqrt{2\pi}}\right)\beta+\frac{1}{2\pi}}+\frac{1}{\sqrt{2\pi}}\right)^{2}.\label{eq:b2}
\end{align}
This is because, 
\begin{align}
\mathbf{W}_{1}[f] & =\mathbb{E}\Big[\sum_{i=1}^{n}\hat{f}_{i}X_{i}f(\mathbf{X})\Big]\label{eq:-45-1}\\
 & \leq\mathbb{E}\Big[\sum_{i=1}^{n}\hat{f}_{i}X_{i}\,1\{\sum_{i=1}^{n}\hat{f}_{i}X_{i}\geq0\}\Big]\label{eq:-43}\\
 & =\frac{1}{2}\mathbb{E}\Big[|\sum_{i=1}^{n}\hat{f}_{i}X_{i}|\Big]\nonumber \\
 & \leq\frac{1}{2}\Big(\sqrt{\frac{2}{\pi}}\sqrt{\mathbf{W}_{1}[f]}+(1-\sqrt{\frac{2}{\pi}})\beta\Big),\label{eq:-31}
\end{align}
where \eqref{eq:-43} follows since if we relax $(\hat{f}_{i})_{i\in[n]}$
and $f$ to be independent quantities, then given $(\hat{f}_{i})_{i\in[n]}$,
the Boolean function $f:\mathbf{x}\mapsto1\{\sum_{i=1}^{n}\hat{f}_{i}x_{i}\geq0\}$
maximizes the expectation in \eqref{eq:-45-1} over all Boolean functions,
and \eqref{eq:-31} follows from the following variant of Khintchine's
inequality proven by König, Schütt, and Tomczak-Jaegermann \cite{konig1999projection}:
\[
\Big|\mathbb{E}\Big[\Big|\sum_{i=1}^{n}c_{i}X_{i}\Big|\Big]-\sqrt{\frac{2}{\pi}}\left\Vert \mathbf{c}\right\Vert _{2}\Big|\leq\Big(1-\sqrt{\frac{2}{\pi}}\Big)\left\Vert \mathbf{c}\right\Vert _{\infty}
\]
with $\mathbf{c}:=(c_{1},c_{2},\ldots,c_{n})$. Solving the inequality
in \eqref{eq:-31}, we obtain \eqref{eq:b2}. 

Combining \eqref{eq:b1} and \eqref{eq:b2} yields $W(a)\le\tilde{\chi}(a)$.
\end{IEEEproof}

Numerical results verify that for $a\in[0.42,0.5]$, 
\[
\tilde{\chi}(a)=\max_{\beta\in[0.4,a]}\frac{1}{4}\left(\sqrt{\chi(a+\beta)}+\sqrt{\chi(a-\beta)}\right)^{2}+\beta^{2}.
\]
Denote $\beta^{*}(a)$ as the maximizer for this optimization. Note
that $\beta^{*}(a)<a$ since the derivative of the objective function
above with respect to $\beta$ goes to $-\infty$ as $\beta\uparrow a$.
But it is numerically verified that $\epsilon^{*}(a):=a-\beta^{*}(a)\downarrow0$
as $a\uparrow1/2$. 

Based on these numerical results, the bound in \eqref{eq:-19} is
asymptotically tight as $a\uparrow1/2$, due to the following arguments.
Denote $f(\mathbf{x})=\frac{1+x_{1}}{2}1_{A^{c}}(x_{2},\ldots,x_{n})+\frac{1-x_{1}}{2}1_{B}(x_{2},\ldots,x_{n})$
with some $A\subseteq\{-1,1\}^{n-1}$ and $B\subseteq\{-1,1\}^{n-1}$
respectively attaining  $W^{(n-1)}(1-a-\beta^{*}(a))$ and $W^{(n-1)}(a-\beta^{*}(a))$.
Then, it holds that 
\[
\mathbf{W}_{1}[f]=\frac{1}{4}\left(\sqrt{W^{(n-1)}(1-a-\beta^{*}(a))}+\sqrt{W^{(n-1)}(a-\beta^{*}(a))}\right)^{2}+\beta^{*}(a)^{2}.
\]
Given $a\le1/2$, letting $n\to\infty$, it holds that 
\[
W(a)\ge\frac{1}{4}\left(\sqrt{W(1-a-\beta^{*}(a))}+\sqrt{W(a-\beta^{*}(a))}\right)^{2}+\beta^{*}(a)^{2},
\]
where $W(1-a-\beta^{*}(a))\sim\chi(1-a-\beta^{*}(a))=\chi(a+\beta^{*}(a))$
and $W(a-\beta^{*}(a))\sim\chi(a-\beta^{*}(a))$ as  $a\uparrow1/2$.

\subsection{\label{subsec:Comparisons}Comparisons }

 An existing bound on $\mathbf{W}_{1}$ proven by the linear programming
method \cite{yu2019improved,fu2001minimum} is 
\begin{equation}
W(a)\le\varphi_{\mathrm{LP}}(a):=\begin{cases}
2a^{2}(\frac{1}{\sqrt{a}}-1), & 0<a\leq\frac{1}{4}\\
\frac{a}{2}, & \frac{1}{4}\leq a\leq\frac{1}{2}
\end{cases}.\label{eq:-17}
\end{equation}
This bound can be improved by using the Khintchine inequality, which
was shown by an anonymous reviewer of the paper \cite{yu2019improved}.
 Our bound in Theorem~\ref{thm:bound} can be further improved by
rechoosing $\varphi$ to $\varphi(t)=\varphi_{\mathrm{LP}}(\tilde{t})$
or the bound proven by the Khintchine inequality, where $\tilde{t}:=\min\{t,1-t\}$.
 Numerical results show that this new bound is only slightly better
than the bound in Theorem~\ref{thm:bound} when $a$ is small. Besides
the bounds mentioned above, there are another two bounds better than
Chang's bound in certain regions: a bound derived by hypercontractivity
given in \cite{yu2021non} and a bound derived by probabilistic inequalities
given in \cite{heilman2025fourier}. 

We now focus on the case $a=1/8$ to compare all the bounds. For this
case, the degree-$1$ Fourier weight of the ($n-3)$-subcube is $\mathbf{W}_{1}[C_{n-3}]=3/64=0.046875$,
which is larger than $J(1/8)$, the degree-$1$ Fourier weight of
Hamming balls of size $1/8$. In contrast, for $a=1/8$, our bound
yields that $\chi(1/8)=0.0505062\ldots$ Our bound is $7.74652\%$ larger
than $\mathbf{W}_{1}[C_{n-3}]$. The bound in \eqref{eq:-17} yields
$\varphi_{\mathrm{LP}}(1/8)=0.0571383\ldots$ The improved version
of our bound mentioned above evaluated at $a=1/8$ is $0.0495142\ldots$
which is $5.63033\%$ larger than $\mathbf{W}_{1}[C_{n-3}]$, better
than other bounds.  

We now compare our bound with Chang's bound in the asymptotic setting
as $a\to0$.  Although Theorem \ref{thm:small-density} shows the optimality of Hamming balls in this regime, there is no explicit bound provided, more precisely, the threshold $a_0$ is required to be sufficiently small, which makes the theorem less applicable to other problems. Thus, we now compare existing bounds with the
explicit bounds given in Theorem~\ref{thm:bound} and \ref{thm:bound2}. 
For comparison, we let $a=e^{-t}$ and consider the
asymptotics of bounds as $t\to\infty$. Denote Chang's bound as $\varphi_{\textrm{Chang}}(a):=2a^{2}\ln\frac{1}{a}$.
For this case, 
\[
\ln\varphi_{\textrm{Chang}}(e^{-t})=-2t+\ln(2t).
\]
By Taylor's expansion, 
\begin{align*}
\ln J(e^{-t}) & =-2t+\ln(2t)-\frac{\ln(2\pi)}{2t}+O\left(\frac{1}{t^{2}}\right),\\
\ln\chi(e^{-t}) & =-2t+\ln\left(2t\right)-\frac{\ln\left(\frac{2\pi}{w^{2}}\right)}{2t}+O\left(\frac{1}{t^{2}}\right).
\end{align*}
We can see that Chang's bound coincides to the first two order terms
of $\ln J(e^{-t})$, while our bound also provides a correct order
for the third order term of $\ln J(e^{-t})$ although the factors
do not coincide. 

\begin{figure}
\centering \includegraphics[width=0.7\columnwidth]{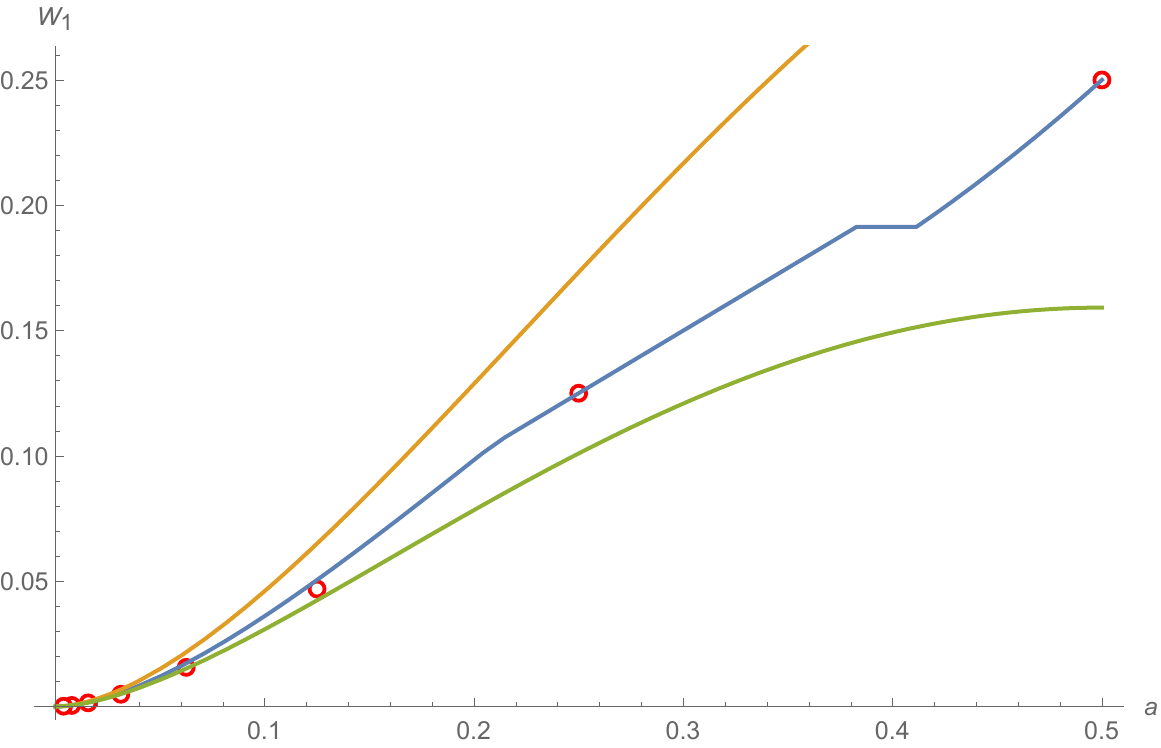}

\includegraphics[width=0.7\columnwidth]{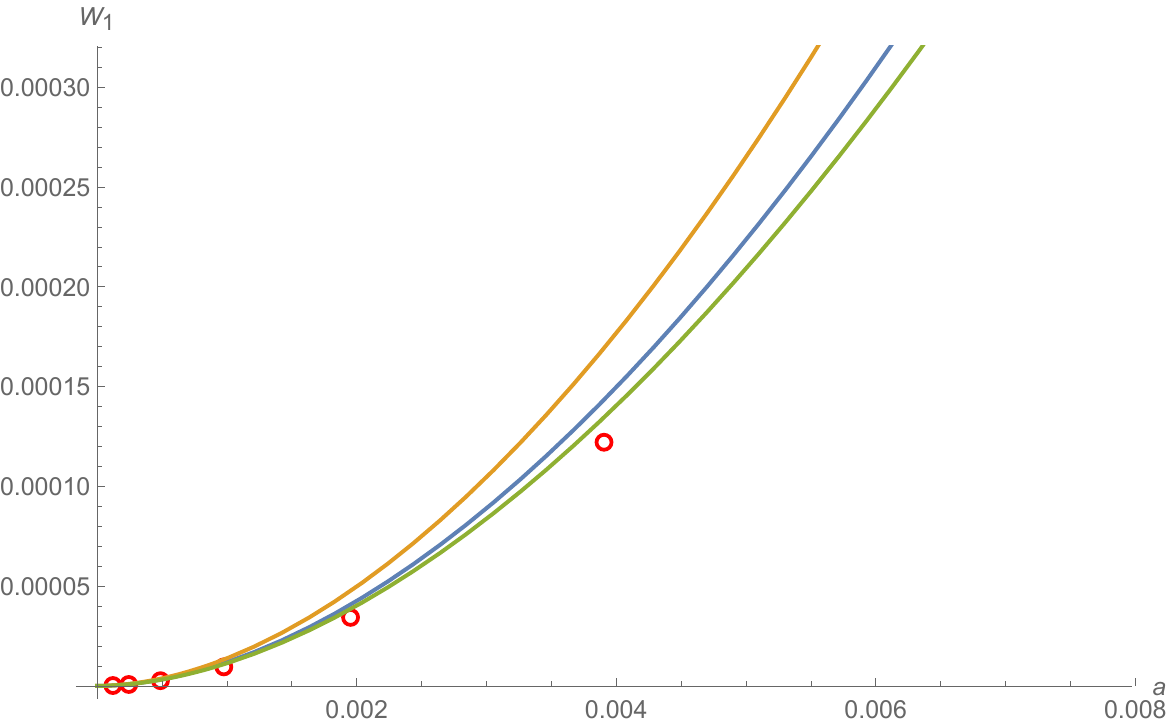}

\caption{\label{fig:bounds}Comparisons of Chang's bound in Theorem~\ref{thm:level-1-bound},
the bounds in Theorem~\ref{thm:bound} and \ref{thm:bound2}, $\mathbf{W}_{1}$
of subcubes, and $\mathbf{W}_{1}$ of Hamming balls. In both the figures,
the solid curves from top to bottom are respectively Chang's bound
in Theorem~\ref{thm:level-1-bound}, the minimum of the bounds in
Theorem~\ref{thm:bound} and \ref{thm:bound2}, and $\mathbf{W}_{1}$
of Hamming balls. The red circles are $\mathbf{W}_{1}$ of subcubes.}
\end{figure}

We compare Chang's bound in Theorem~\ref{thm:level-1-bound}, the
bounds in Theorem~\ref{thm:bound} and \ref{thm:bound2}, $\mathbf{W}_{1}$
of subcubes, and $\mathbf{W}_{1}$ of Hamming balls in Fig. \ref{fig:bounds}.
From the figure, we can see that our bounds are very close to the
lower bound---the maximum of $\mathbf{W}_{1}$ of subcubes and $\mathbf{W}_{1}$
of Hamming balls, no matter for small or large $a$. 

\section{Application to the Friedgut--Kalai--Naor Theorem}\label{sec:fkn}

The Friedgut--Kalai--Naor (FKN) theorem concerns about which Boolean
functions $f$ on the discrete cube have Fourier coefficients concentrated
at the lowest two levels. It states that such Boolean functions are
close to either a constant function (i.e., $f=0$ or $1$) or a dictator
function ($f=1\{x_{i}=1\}$ or $1\{x_{i}=-1\}$). Here we aim at improving
the FKN theorem by focusing on the class of Boolean functions with
a given mean. For this case, the FKN theorem can be also formulated
as a theorem about maximizing the first-order Fourier weight of a
Boolean function given the maximum of its correlations to all dictator
functions. We next provide the formulation.

For $a,\beta\in2^{-n}[2^{n}]$, define

\begin{equation}
W^{(n)}(a,\beta):=\max_{f:\{-1,1\}^{n}\to\{0,1\}:\mathbb{E}f=a,\max_{i\in[n]}|\hat{f}_{i}|=\beta}\mathbf{W}_{1}[f].\label{eq:W-3}
\end{equation}
Define $W(a,\beta):=\sup_{n\ge1}W^{(n)}(a,\beta)$. Here without ambiguity,
we use the same notation $W^{(n)}$ but with different numbers of
parameters to denote two different functions $W^{(n)}(a,\beta)$ and
$W^{(n)}(a)$. The quantity $W^{(n)}(a,\beta)$ was implicitly studied
by Friedgut, Kalai, and Naor \cite{friedgut2002boolean} who showed
that for $a=\frac{1}{2}$, $W^{(n)}(\frac{1}{2},\beta)\to\frac{1}{4}$
if and only if $\beta\to\frac{1}{2}$.
\begin{lem}
It holds that for $\beta\in[0,1/2]$, 
\[
W^{(n)}(\frac{1}{2},\beta)\le\beta^{2}+W^{(n-1)}(\frac{1}{2}-\beta),
\]
and the equality above holds for $\beta\ge1/4$. Moreover,  $W^{(n)}(\frac{1}{2},\beta)$
for $\beta\ge1/4$ is attained by 
\[
f(\mathbf{x})=\frac{1+x_{1}}{2}1_{A^{c}}(x_{2},\ldots,x_{n})+\frac{1-x_{1}}{2}1_{A}(x_{2},\ldots,x_{n})
\]
for some $A\subseteq\{-1,1\}^{n-1}$ attaining  $W^{(n-1)}(\frac{1}{2}-\beta)$. 
\end{lem}
\begin{IEEEproof}
This is a consequence of Lemma 2 in \cite{yu2023phi}. 
\end{IEEEproof}
By this lemma, determining the function $\beta\in[1/4,1/2]\mapsto W(\frac{1}{2},\beta)$
is equivalent to determining the function $a\in[0,1/4]\mapsto W(a)$.
Combining this lemma with Theorem~\ref{thm:bound} yields the following
bound. 
\begin{thm}
It holds that $W(\frac{1}{2},\beta)\le\beta^{2}+\chi(\frac{1}{2}-\beta)$,
where the function $\chi$ is given in \eqref{eq:J} with $\varphi(t)=t/2$
and $T=0.21$. 
\end{thm}

This is an improvement of the existing bound in \cite{yu2023phi}
derived by using the level-1 weight bound in Theorem~\ref{thm:level-1-bound}.
When $\beta\to1/2$, the bound in this paper is asymptotically sharp
in certain sense. When $\beta<1/4$, the bound derived in \cite{yu2023phi}
given below is better: 
\begin{equation}
W(\frac{1}{2},\beta)\leq\frac{1}{4}\Big(\sqrt{4(\frac{1}{2}-\frac{1}{\sqrt{2\pi}})\beta+\frac{1}{2\pi}}+\frac{1}{\sqrt{2\pi}}\Big)^{2}.\label{eq:-21-2}
\end{equation}

\section{Average Distance in Euclidean Spaces}\label{sec:euclidean}

Let $P$ be a probability measure on $\mathbb{R}^{n}$. Let $\mathbf{X},\mathbf{Y}$
be independent random vectors both obeying $P$. Define the average
distance of $A$ w.r.t. $P$ as 
\[
\mathbf{D}(A):=\sqrt{\mathbb{E}\Big[\sum_{i=1}^{n}(X_{i}-Y_{i})^{2}\Big|\mathbf{X}\in A,\mathbf{Y}\in A\Big]}.
\]
Then, it is easy to see that 
\[
\mathbf{D}^{2}(A)=2\sum_{i=1}^{n}\operatorname{Var}(X_{i}\Big|\mathbf{X}\in A).
\]
So, minimizing the average distance of a set is equivalent to minimizing
the conditional variance given the set. 
\begin{thm}
Euclidean balls minimize the average distance over all sets with the
same measure. 
\end{thm}
\begin{IEEEproof}
The minimizers must be a Euclidean ball, not only for Gaussian measures
but also for any other measures. This is because, for any set $A$,
the barycenter is denoted by $\mathbf{c}=(\mathbb{E}[X_{i}\Big|\mathbf{X}\in A])_{i\in[n]}=(\mathbb{E}_{Q}[X_{i}])_{i\in[n]}$,
where $Q:=P(\cdot|A)$. Then, 
\[
\mathbf{D}^{2}(A)=2\mathbb{E}_{Q}\Big[\sum_{i=1}^{n}(X_{i}-c_{i})^{2}\Big].
\]
Denote $B_{r}(\mathbf{c})$ as the ball of radius $r$ at center $\mathbf{c}$
having probability equal to $P(A)$. Then, if move all points from
$A\backslash B_{r}(\mathbf{c})$ into $B_{r}(\mathbf{c})$, then the
value of $\mathbb{E}_{Q}\left[\sum_{i=1}^{n}(X_{i}-c_{i})^{2}\right]$
reduces. So, the minimizers must be a Euclidean ball or a set having
symmetric difference with a Euclidean ball equal to zero. 
\end{IEEEproof}
The above theorem implies Theorem~\ref{thm:optimality}, i.e., linear
threshold sets minimize the average distance in Boolean hypercube
(equivalently, maximize the degree-1 Fourier weight). This can be
seen from that when $P$ is the uniform distribution over $\{\pm1\}^{n}$,
it holds that for any $A\subseteq\{\pm1\}^{n}$, 
\begin{align*}
\mathbf{D}^{2}(A) & =\mathbb{E}\Big[\sum_{i=1}^{n}(X_{i}-Y_{i})^{2}\Big|\mathbf{X}\in A,\mathbf{Y}\in A\Big]\\
 & =4\mathbb{E}\Big[d_{\mathrm{H}}(\mathbf{X},\mathbf{Y})\Big|\mathbf{X}\in A,\mathbf{Y}\in A\Big],
\end{align*}
where $d_{\mathrm{H}}$ is the Hamming distance. 

For the Lebesgue measure on $\mathbb{R}^{n}$, the sets minimizing
the average distance are all Euclidean balls of given size.  For
the standard Gaussian measure on $\mathbb{R}^{n}$, the set minimizing
the average distance is unique and is the ball at $\mathbf{0}$ of
a given measure. This is because, for any ball $B_{r}(\mathbf{c})$,
\[
\mathbf{D}^{2}(B_{r}(\mathbf{c}))=2\mathbb{E}_{Q}\left[Y\right]=2\int_{0}^{\infty}(1-F(y))\mathrm{d}y,
\]
where $Y:=\sum_{i=1}^{n}(X_{i}-c_{i})^{2}$, and $F$ is the conditional
CDF of $Y$ under the condition $\mathbf{X}\in B_{r}(\mathbf{c})$,
i.e., $F(y):=P\{Y\le y|\mathbf{X}\in B_{r}(\mathbf{c})\}=P\{Y\le y\}/P(B_{r}(\mathbf{c}))$
for $y\le r$. So, a set having larger $F$ would have smaller $\mathbf{D}(A)$.
Let $t$ be chosen such that $P(B_{t}(\mathbf{0}))=B_{r}(\mathbf{c})$.
For the set $B_{t}(\mathbf{0})$, the conditional CDF turns into $G(y)=P\{\sum_{i=1}^{n}X_{i}^{2}\le y\}/P(B_{r}(\mathbf{0}))$.
For the standard Gaussian measure $P$, it holds that $F(y)\le G(y),\forall y$,
since 
\begin{align*}
P\{\sum_{i=1}^{n}(X_{i}-c_{i})^{2}\le y\} & \le P\{\sum_{i=1}^{n}X_{i}^{2}\le y\}.
\end{align*}
 So, the minimizer is $B_{t}(\mathbf{0})$. 

\section{Chang's Original Lemma for $\mathbb{F}_{2}^{n}$}\label{sec:chang}

Theorem~\ref{thm:level-1-bound} is in fact not the original version
of Chang's lemma for $\mathbb{F}_{2}^{n}$. The original version of
Chang's lemma focuses on the estimation of the dimension of the space
spanned by large Fourier coefficients, which was widely used in additive
combinatorics. In this section, we also consider the space $\mathbb{F}_{2}^{n}$,
the $n$-th power of the finite field of size $2$. The spaces $\mathbb{F}_{2}^{n}$
and $\{\pm1\}^{n}$ admit a natural bijection $y_{i}=(-1)^{x_{i}},i\in[n]$
for $\mathbf{x}\in\mathbb{F}_{2}^{n}$ and $\mathbf{y}\in\{\pm1\}^{n}$.
So, we do not differentiate Boolean functions $f:\mathbb{F}_{2}^{n}\to\mathbb{F}_{2}$
and $f:\{\pm1\}^{n}\to\{0,1\}$. We also use $\hat{f}(\mathbf{y}),\mathbf{y}\in\mathbb{F}_{2}^{n}$
to denote Fourier coefficients of $f$.

Suppose that $f:\mathbb{F}_{2}^{n}\to\mathbb{F}_{2}$ is a Boolean
function, and $\epsilon\in(0,1]$. Then we define the \emph{$\epsilon$-spectrum}
of $f$ as 
\[
\mathrm{Spec}_{\epsilon}(f):=\{\mathbf{y}\in\mathbb{F}_{2}^{n}:|\hat{f}(\mathbf{y})|>\epsilon\|f\|_{1}\}.
\]
We now estimate the dimension of the space spanned by $\mathrm{Spec}_{\epsilon}(f)$,
i.e., $\mathrm{dim}(\mathrm{Span}(\mathrm{Spec}_{\epsilon}(f)))$. 

Let $A$ be the support of $f$. Let $\gamma_{1},\cdots,\gamma_{d}$
be a maximal set of linearly independent vectors in $\mathrm{Spec}_{\epsilon}(f)$.
That is, $d=\mathrm{dim}(\mathrm{Span}(\mathrm{Spec}_{\epsilon}(f)))$.
Let $M\in\mathbb{F}_{2}^{n\times n}$ be an invertible matrix (a matrix
consisting of $n$ linearly independent vectors) whose first $d$
rows are $\gamma_{1},\cdots,\gamma_{d}$. Let $B=MA=\{M\mathbf{x}:\mathbf{x}\in A\}$
and $g=1_{B}$. Then, it is easy to see that $\mathbb{E}f=\mathbb{E}g$,
and $\hat{f}(\gamma_{1})=\hat{g}(e_{1}),\cdots,\hat{f}(\gamma_{d})=\hat{g}(e_{d})$.
Here $e_{i}=(0,\cdots,0,1,0,\cdots,0)$, with the $1$ in the $i$-th
position. This observation follows since for $1\le i\le d$, 
\begin{align*}
\hat{g}(e_{i}) & =\frac{1}{q^{n}}\sum_{\mathbf{x}\in B}(-1)^{\langle\mathbf{x},e_{i}\rangle}=\frac{1}{q^{n}}\sum_{\mathbf{x}\in A}(-1)^{\langle M\mathbf{x},e_{i}\rangle}=\frac{1}{q^{n}}\sum_{\mathbf{x}\in A}(-1)^{\langle\mathbf{x},M^{\top}e_{i}\rangle}\\
 & =\frac{1}{q^{n}}\sum_{\mathbf{x}\in A}(-1)^{\langle\mathbf{x},\gamma_{i}\rangle}=\hat{f}(\gamma_{i}).
\end{align*}
So, estimating the dimension of the space spanned by $\mathrm{Spec}_{\epsilon}(f)$
is equivalent to estimating the number of large level-1 Fourier coefficients.
In fact, the quantity $W(a)$ can be used to bound this dimension. 
\begin{lem}
\label{lem:It-holds-that}Let $f:\mathbb{F}_{2}^{n}\to\mathbb{F}_{2}$
have expectation $a$. Then, 
\[
\mathrm{dim}(\mathrm{Span}(\mathrm{Spec}_{\epsilon}(f)))\le\frac{W(a)}{a^{2}\epsilon^{2}}.
\]
\end{lem}
\begin{IEEEproof}
 Let $f:\mathbb{F}_{2}^{n}\to\mathbb{F}_{2}$ be a function such
that $\mathbb{E}f=a,\,|\hat{f}_{i}|>a\epsilon,\forall i\in[d]$. Then,
$W(a)\ge\sum_{i\in[d]}|\hat{f}_{i}|^{2}>da^{2}\epsilon^{2}.$ 
\end{IEEEproof}
The following theorem is the original lemma derived by Chang, which
is important in additive combinatorics. This theorem follows by combining
Lemma \ref{lem:It-holds-that} and Theorem~\ref{thm:level-1-bound}. 
\begin{thm}[Chang's Lemma]
\label{thm:chang} \cite{chang2002polynomial} Let $f:\mathbb{F}_{2}^{n}\to\mathbb{F}_{2}$
have expectation $a$. Then, 
\begin{equation}
\mathrm{dim}(\mathrm{Span}(\mathrm{Spec}_{\epsilon}(f)))\le2\epsilon^{-2}\log\frac{1}{a}.\label{eq:chang}
\end{equation}
\end{thm}
Although the bounds in Theorems \ref{thm:bound} and \ref{thm:bound2}
can be applied in the same way to strengthen Chang's lemma, in the
following, we prove a sharp version of Chang's lemma. To this end,
we reformula the problem in a different way. Given $a$ and $k$,
we ask what is the maximum possible value of $\epsilon$ such that
$a=\mathbb{E}f$ and $k=\mathrm{dim}(\mathrm{Span}(\mathrm{Spec}_{\epsilon}(f)))$
for some Boolean function $f$.
\begin{thm}
\label{thm:sharpchang}Let $\epsilon>0$. Let $f:\mathbb{F}_{2}^{n}\to\mathbb{F}_{2}$
have expectation $a=\mathbb{E}f$ and $k=\mathrm{dim}(\mathrm{Span}(\mathrm{Spec}_{\epsilon}(f)))$.
Let $h(\mathbf{x})=1\{\sum_{i=1}^{k}x_{i}>b\}+\lambda\cdot1\{\sum_{i=1}^{k}x_{i}=b\}$
for some $b\in\mathbb{N}$ and $\lambda\in(0,1]$ such that $\mathbb{E}h=a$.
 Then, $\hat{h}_{i},i\in[k]$ are all equal, and $\epsilon\le a^{-1}\hat{h}_{1}$. 
\end{thm}

\begin{IEEEproof}
Let $A$ be the support of $f$. Without loss of generality, we assume
$\hat{f}_{i}\ge0$ for all $i\in[n]$, since otherwise, we can change
the sign of the $i$-th coordinate for all points in $A$. Observe
that $\hat{f}_{i}=\mathbb{E}[f(\mathbf{X})X_{i}]=\frac{a}{|A|}\sum_{\mathbf{x}\in A}x_{i}$.
So, $a^{-1}(\hat{f}_{i})_{i\in[n]}$ is the barycenter of $A$. Let
$\beta=\min_{i\in[k]}\hat{f}_{i}$. Then, we claim that there is a
set $B$ such that $\mu(B)=\mu(A)$ and $\hat{g}_{i}=\beta,i\in[k]$,
where $g=1_{B}$. We now prove this claim. 

For $i\in[k]$, denote $A_{1}=A\cap\{x_{i}=1\}$ and $A_{-1}=A\cap\{x_{i}=-1\}$,
which respectively have sizes $2^{n-1}(a+\hat{f}_{i})$ and $2^{n-1}(a-\hat{f}_{i})$.
We choose arbitrary $2^{n-1}(\hat{f}_{i}-\beta)$ points from $A_{1}\backslash A_{-1}$,
and then change the $i$-th coordinates of them from $1$ to $-1$.
 For each $i\in[k]$, we perform this change once. We will obtain
the desired set $B$. 

The barycenter of $B$ is $\mathbf{w}=(\beta/a,...,\beta/a,w_{k+1},...,w_{n})$
for some $w_{k+1},...,w_{n}$. We continue moving points from $B\cap\{\sum_{i=1}^{k}x_{i}<b\}$
to $\{\sum_{i=1}^{k}x_{i}>b\}$ until the set $\{\sum_{i=1}^{k}x_{i}>b\}$
is filled up; we then uniformly allocate all the total mass of the
remaining points in $B\cap\{\sum_{i=1}^{k}x_{i}<b\}$ to all points
in $\{\sum_{i=1}^{k}x_{i}=b\}$. So, the  ``indicator'' of the
resultant set is just the function $h$. During the whole procedure
of this construction, the component of the barycenter of the set along
the direction $(1,...,1,0,...,0)$ (consisting of $k$ ones and $n-k$
zeros) is always increasing. So, it holds that $\sum_{i=1}^{k}\hat{h}_{i}^{2}\ge\sum_{i=1}^{k}\hat{g}_{i}^{2}$.
   Since $\hat{g}_{i},i\in[k]$ are all equal to $\beta$ and $\hat{h}_{i},i\in[k]$
are all equal, it holds that $\beta\le\min_{i\in[k]}\hat{h}_{i}$.
\end{IEEEproof}
As a consequence of the theorem above, we have the following corollary. 
\begin{cor}
\label{cor:Under-the-condition}Under the condition same as the one
in Theorem~\ref{thm:sharpchang}. Suppose that there is a set $C\subseteq\{\pm1\}^{n}$
such that $\mu(C)=a$, $\{\sum_{i=1}^{k}x_{i}>b\}\subseteq C\subseteq\{\sum_{i=1}^{k}x_{i}\ge b\}$,
and $\hat{g}_{i},i\in[k]$ are all equal where $g=1_{C}$. Then, $\hat{g}_{1}=\hat{h}_{1}$
and $\epsilon\le a^{-1}\hat{g}_{1}$.   In particular, when $a=2^{-k}{k \choose \le r}$
for some $r$, the set $C$ can be chosen as a Hamming ball $\{\sum_{i=1}^{k}x_{i}\ge k-2r\}$,
and for this case, $\epsilon\le\frac{{k-1 \choose r}}{{k \choose \le r}}$,
where ${k \choose \le r}={k \choose 1}+...+{k \choose r}$. 
\end{cor}
\begin{IEEEproof}
It is easy to see that $\hat{p}_{i},\hat{q}_{i},i\in[k]$ are all
equal where $p=1_{E}$ and $q=\lambda\cdot1\{\sum_{i=1}^{k}x_{i}=b\}$
with $E=C\cap\{\sum_{i=1}^{k}x_{i}=b\}$ and $\lambda=\mu(E)$. This
is because, by assumption, $\hat{p}_{i},i\in[k]$ are all equal, and
hence, $\hat{p}_{i}=\frac{1}{k}\sum_{i\in[k]}\hat{p}_{i}=2^{-n}\sum_{\mathbf{x}\in E}\frac{1}{k}\sum_{i=1}^{k}x_{i}=2^{-n}\sum_{\mathbf{x}\in E}\frac{b}{k}=\frac{\lambda b}{k}=\hat{q}_{j}$
for any $i,j\in[k]$. 
\end{IEEEproof}
Corollary \ref{cor:Under-the-condition} implies that Hamming balls
have the largest dimension of the space spanned by their large Fourier
coefficients. So, for fixed $a$, we have $\epsilon\lesssim\frac{I(a)}{a\sqrt{k}}$
as $k\to\infty$. In other words, $k\lesssim\frac{J(a)}{a^{2}\epsilon^{2}}$
as $\epsilon\to0$.  Here, we denote $f(x)\lesssim g(x)$ as $x\to\infty$
if $\limsup_{x\to\infty}\frac{f(x)}{g(x)}\le1$. 

Chang's bound was previously improved to a strengthened but non-sharp
version by Shkredov \cite{shkredov2008sets} for $\mathbb{Z}_{N}$.
He then applied his new bound to additive combinatorics, strengthening
some result of Chang \cite{chang2002polynomial}. By replacing our
sharp version of Chang's bound with his, one can obtain a further
strengthening of Chang's result for $\mathbb{F}_{2}^{n}$. Chang's
bound was also improved in \cite{chakraborty2020tight}. Our sharp
version of Chang's bound can be easily used to slightly refine Bogolyubov\textquoteright s
lemma and Green's lemma \cite[Theorem 20]{green2003restriction};
see \cite{sandersanalysis} for details. Besides, Chang's lemma has
also found numerous applications in other fields; see a comprehensive
review in \cite{chakraborty2020tight}.

\section{Open Problems}\label{sec:open-problems}

Problem~\ref{prob:W} is resolved by Theorem~\ref{thm:small-density}.
Two further questions concerning the Fourier weight remain.
\begin{problem}
What is the value of $W(1/8)$?
\end{problem}
As discussed in Section~\ref{subsec:Comparisons}, the best bound obtained
above is $5.63033\%$ larger than $\mathbf{W}_{1}[C_{n-3}]$. 

It is known that, as $\beta\to0$, $W(1/2,\beta)\to\frac{1}{2\pi}$.
It is therefore natural to ask for the rate of this convergence. 
\begin{problem}
Is it true that $W(1/2,\beta)=\frac{1}{2\pi}+O(\beta^{2})$ as $\beta\to0$?
\end{problem}
This conclusion holds if $W(1/2,\beta)$ is attained by the LTF $f(\mathbf{x})=1\{\sum_{i=1}^{k}x_{i}\ge0\}$
with odd $k$ chosen so that  $|\hat{f}_{i}|\approx\beta,1\le i\le k$. 

\textit{Updated Note:} Problem~3 has since been resolved affirmatively by the author in \cite{yu2026sharp}.

\appendices
\section{\label{sec:-Proof-of}Alternative Proof of Theorem~\ref{thm:optimality}}

Define for $\mathbf{x}\in\{\pm1\}^{n}$, $g(\mathbf{x})=\sum_{i=1}^{n}\hat{f}_{i}x_{i}.$
Observe 
\begin{align*}
\mathbf{W}_{1}[f]=\sum_{i=1}^{n}\hat{f}_{i}^{2} & =\mathbb{E}[g(\mathbf{X})f(\mathbf{X})]\le\frac{1}{2^{n}}\sum_{\mathbf{x}}g(\mathbf{x})\ell(\mathbf{x})\\
 & =\sum_{i=1}^{n}\hat{f}_{i}\hat{\ell}_{i}\leq\sqrt{(\sum_{i=1}^{n}\hat{f}_{i}^{2})(\sum_{i=1}^{n}\hat{\ell}_{i}^{2})}\\
 & =\sqrt{\mathbf{W}_{1}[f]\mathbf{W}_{1}[\ell]},
\end{align*}
where the first inequality follows since given $g$, a Boolean-valued
function $\ell$ of expectation $a$ maximizes $\sum_{x^{n}}g(x^{n})\ell(x^{n})$
if and only if $f_{>}\le\ell\le f_{\ge}$ and $\mathbb{E}\ell=a$.
Hence $\mathbf{W}_{1}[f]\le\mathbf{W}_{1}[\ell]$. By the optimality
of $f$, equality holds, and therefore $\ell$ is also optimal. Hence,
any Boolean-valued function $\ell$ such that $f_{>}\le\ell\le f_{\ge}$
and $\mathbb{E}\ell=a$ is optimal.

The optimality of $f$ and $\ell$ further implies that every inequality
in the chain above is an equality. Hence, $f_{>}\le f\le f_{\ge}$.
Hence, the set of optimal Boolean-valued functions is $\{\textrm{Bool }\ell:f_{>}\le\ell\le f_{\ge},\mathbb{E}\ell=a\}$.
Consequently, the definitions of $f_{>}$ and $f_{\ge}$ do not depend on a
particular optimizer $f$; otherwise the resulting class of optimizers would
not be uniquely determined.

\section{\label{sec:Proof-of-Lemma}Proof of Lemma \ref{lem:Gamma}}

 Observe that 
\begin{align*}
\lim_{a_{1}\downarrow T}\Gamma_{2}'(a_{0},a_{1}) & =\frac{1}{2}\left(I_w(a_{0})+\sqrt{\frac{T}{2}}\right)\frac{1}{2\sqrt{2T}}+\frac{1}{2}(T-a_{0})-I_w\left(\frac{a_{0}+T}{2}\right)I_w'\left(\frac{a_{0}+T}{2}\right)\\
 & \le\frac{1}{2}\left(I_w(a_{0})+I_w(T)\right)I_w'(T)+\frac{1}{2}(T-a_{0})-I_w\left(\frac{a_{0}+T}{2}\right)I_w'\left(\frac{a_{0}+T}{2}\right)\\
 & =\lim_{a_{1}\uparrow T}\Gamma_{2}'(a_{0},a_{1})\le0,
\end{align*}
where the first inequality follows since $I_w(T)=\sqrt{\frac{T}{2}}$
and $\frac{1}{2\sqrt{2T}}=0.771517\ldots\le1.02231\ldots=I_w'(T)$, and the
second inequality follows from \eqref{eq:-11}. 

We next prove $\Gamma_{2}'(a_{0},2T-a_{0})\le0$ for $0.02\le a_{0}\le T$.
Observe that 
\begin{align*}
\Gamma_{2}'(a_{0},2T-a_{0}) & =\frac{1}{2}\left(I_w(a_{0})+\sqrt{\frac{2T-a_{0}}{2}}\right)\frac{1}{2\sqrt{2(2T-a_{0})}}+(T-a_{0})-I_w\left(T\right)I_w'\left(T\right)\\
 & =\frac{1}{2}\left(\frac{I_w(a_{0})}{2\sqrt{2(2T-a_{0})}}+\frac{1}{4}\right)+(T-a_{0})-I_w\left(T\right)I_w'\left(T\right)\\
 & =:h(a_{0}),
\end{align*}
and 
\begin{align*}
h'(a_{0}) & =\frac{2I_w'(a_{0})(2T-a_{0})+I_w(a_{0})}{8\sqrt{2}(2T-a_{0})^{3/2}}-1.
\end{align*}
We claim that $h'(a_{0})\le0$ for $0.02\le a_{0}\le T$. To prove
this claim, we only need show that 
\[
\eta(a_{0}):=2I_w'(a_{0})(2T-a_{0})+I_w(a_{0})-8\sqrt{2}(2T-a_{0})^{3/2}\le0.
\]
Observe that 
\begin{align}
\eta(a_{0}) & \le2I_w'(a_{0})(2T-a_{0})+I_w(T)-8\sqrt{2}(2T-a_{0})^{3/2}\nonumber \\
 & =2(2T-a_{0})\left(I_w'(a_{0})-4\sqrt{2(2T-a_{0})}\right)+I_w(T),\label{eq:-13}
\end{align}
where the inequality follows because $I_w$ is increasing on $[0.02,T]$.
Moreover, for $0.02\le a_{0}\le T$, it holds that 
\[
2T-a_{0}\ge T,
\]
and 
\begin{align}
I_w'(a_{0})-4\sqrt{2(2T-a_{0})} & \le\max\{I_w'(0.02)-4\sqrt{2(2T-0.02)},I_w'(T)-4\sqrt{2T}\}\label{eq:-12}\\
 & =-1.39698\ldots\nonumber 
\end{align}
where the inequality in \eqref{eq:-12} follows because both $I_w'(a_{0})$
and $-4\sqrt{2(2T-a_{0})}$ are convex (for the former, see \eqref{eq:-5}).
Substituting these two inequalities and $T=0.21$ into \eqref{eq:-13}
yields that $\eta(a_{0})\le-0.262693\ldots<0.$ So, $h'(a_{0})\le0$
for $0.02\le a_{0}\le T$, which further implies $h(a_{0})\le h(0.02)=-0.00549341\ldots<0$.

Finally, we prove $\Gamma_{2}'(a_{0},0.4)\le0$ for $a_{0}\le0.02$.
For $a_{0}\le0.02$, 
\begin{align*}
\Gamma_{2}'(a_{0},0.4) & =\frac{1}{2}\left(I_w(a_{0})+\sqrt{0.2}\right)\frac{1}{2\sqrt{0.8}}+\frac{1}{2}(0.4-a_{0})-I_w\left(\frac{a_{0}+0.4}{2}\right)I_w'\left(\frac{a_{0}+0.4}{2}\right)\\
 & =:h(a_{0}).
\end{align*}
Moreover,
\begin{align}
h'(a_{0}) & =\frac{1}{4\sqrt{0.8}}I_w'(a_{0})-\frac{1}{2}-\frac{1}{2}\left(I_w'\left(\frac{a_{0}+0.4}{2}\right)\right)^{2}-\frac{1}{2}I_w\left(\frac{a_{0}+0.4}{2}\right)I_w''\left(\frac{a_{0}+0.4}{2}\right)\nonumber \\
 & =\frac{1}{4\sqrt{0.8}}I_w'(a_{0})-\frac{1}{2}\left(I_w'\left(\frac{a_{0}+0.4}{2}\right)\right)^{2}\label{eq:-14}\\
 & \ge\frac{1}{4\sqrt{0.8}}I_w'(0.02)-\frac{1}{2}\left(I_w'\left(0.2\right)\right)^{2}\label{eq:-15}\\
 & =0.0544183\ldots>0,\nonumber 
\end{align}
where \eqref{eq:-14} follows since $I_w(a)I_w''(a)=-1$ for all $a\in(0,1)$,
and \eqref{eq:-15} follows because $I_w'$ is decreasing on $[0,0.02]$
(see \eqref{eq:-16}). So, $h(a_{0})\le h(0.02)=-0.00549341\ldots<0$.

\section*{Acknowledgements }

\emph{Generative-AI use disclosure:} ChatGPT 5.6 was used to assist in developing the proofs presented in Section II. The author independently verified all arguments and mathematical claims. The author assumes full responsibility for the correctness and content of the manuscript.

\emph{Funding:} This work was supported by the National Key Research
and Development Program of China under grant 2023YFA1009604 and the
NSFC under grant 62101286.

\bibliographystyle{unsrt}
\bibliography{ref}

\end{document}